\documentclass[12pt]{article}

\usepackage{amssymb}
\usepackage{amsmath}
\usepackage{amscd}

\numberwithin{equation}{section}

\newtheorem{theorem}{Theorem}[section]
\newtheorem{corollary}{Corollary}[section]

\begin{document}

\title{ASYMPTOTIC EXPANSIONS FOR THE SOJOURN TIME DISTRIBUTION IN THE $M/G/1$-PS QUEUE}

\author{
Qiang Zhen\thanks{
Department of Mathematics, Statistics, and Computer Science,
University of Illinois at Chicago, 851 South Morgan (M/C 249),
Chicago, IL 60607-7045, USA.
{\em Email:} qzhen2@uic.edu.}
\and
and
\and
Charles Knessl\thanks{
Department of Mathematics, Statistics, and Computer Science,
University of Illinois at Chicago, 851 South Morgan (M/C 249),
Chicago, IL 60607-7045, USA.
{\em Email:} knessl@uic.edu.\
\newline\indent\indent{\bf Acknowledgement:} This work was partly supported by NSF grant DMS 05-03745.
}}
\date{January 30, 2008 }
\maketitle

\begin{abstract}

We consider the $M/G/1$ queue with a processor sharing server. We study the conditional sojourn time distribution, conditioned on the customer's service requirement, as well as the unconditional distribution, in various asymptotic limits. These include large time and/or large service request, and heavy traffic, where the arrival rate is only slightly less than the service rate. Our results demonstrate the possible tail behaviors of the unconditional distribution, which was previously known in the cases $G=M$ and $G=D$ (where it is purely exponential). We assume that the service density decays at least exponentially fast. We use various methods for the asymptotic expansion of integrals, such as the Laplace and saddle point methods.

\end{abstract}

\newpage

\section{Introduction}

One of the most interesting service disciplines in queueing theory is that of processor sharing (PS). Here every customer in the system gets an equal fraction of the server or processor, and this has the advantage that shorter jobs get served in less time than, say, under the first-in-first-out (FIFO) discipline.

PS queues were introduced in the 1960's by Kleinrock (see \cite{KLanalysis},\cite{KLtime}) and have been the subject of much research over the past $40^+$ years. In these models one of the main measures of performance is a given (also called tagged) customer's sojourn time distribution, conditioned on that customer's service time. The sojourn time is the time the tagged customer leaves the system after being served, assuming the customer arrives at time zero.

We denote by $\mathbf{V}(x)$ the conditional sojourn time, with
$x$ being the service time. If the tagged customer arrived to an
empty system and no further arrivals occurred in the time interval
$[0,x]$, then $\mathbf{V}(x)=x$. But in general $\mathbf{V}(x)>x$
as the tagged customer must share the server. We denote by $b(x)$
the service time density, by $p(t|x)$ the conditional sojourn time
density, and by $p(t)=\int_0^tb(x)\,p(t|x)dx$ the unconditional
sojourn time density. We note that $p(t|x)$ has in general a
probability mass along $t=x$, but $p(t)$ is generally continuous.

The $M/M/1$-PS queue assumes Poisson arrivals and i.i.d. service times with density $b(x)=\mu\, e^{-\mu x}$. In \cite{COF}, Coffman, Muntz and Trotter obtained an expression for the Laplace transform of $p(t|x)$ (i.e., for $\mathbf{E}[e^{-s\mathbf{V}(x)}]$, where $s$ is the Laplace transform variable). In \cite{MO} Morrison removed the conditioning on $x$ and studied $p(t)$ in the heavy traffic limit, where the Poisson arrival rate $\lambda$ is nearly equal to the service rate $\mu$ (thus $\rho=\lambda/\mu\uparrow 1$). Setting $\epsilon=1-\rho$, in \cite{MO} asymptotic results were obtained for the time scales $t=O(1)$, $t=O(\epsilon^{-1})$ and $t=O(\epsilon^{-3})$. Most the mass is concentrated in the range $t=O(\epsilon^{-1})$, and the asymptotic series involves modified Bessel functions. For moderate traffic intensities with $\rho<1$, the tail behavior of $p(t)$ is given by
\begin{equation}\label{1.1}
p(t)\sim C_0\,t^{-5/6}\,e^{-A_0\,t}\,e^{-B_0\,t^{1/3}},\;\;t\rightarrow\infty
\end{equation}
where $A_0=\mu\,(1-\sqrt{\rho})^2$ and the constants $B_0$ and $C_0$ are given in \cite{FL}-\cite{PO}. The result in (\ref{1.1}) was obtained for the $M/M/1$ queue under a random order of service (ROS) discipline, but there is a close connection between the waiting time distribution in the ROS model and the sojourn time distribution in the PS model. This relation, along with some extensions, is explored in \cite{COH} and \cite{BO}. In \cite{ZH} we studied the conditional density $p(t|x)$ for various asymptotic ranges of $x$ and $t$, for both a fixed $\rho<1$ and in the heavy traffic limit where $\rho\uparrow 1$.

A more difficult model is the $M/G/1$-PS queue, where the service density is general. This was analyzed by Yashkov in \cite{YAproc}, \cite{YAmath} and by Ott \cite{OT}. These authors obtained an explicit, albeit complicated, expression for $\mathbf{E}[e^{-s\mathbf{V}(x)}]$. Inverting the Laplace transform leads to an expression for $p(t|x)$ as a contour integral (see (\ref{Ee-sV2})), but the integrand is a nonlinear function of another contour integral, which is in turn defined in terms of the Laplace transform of the service density. In the case of deterministic service times, where $b(x)=\delta(x-1/\mu)$, much more explicit results are available (see \cite{OT} and \cite{EG}). Also, the tail behavior of the unconditional sojourn time density was derived by Egorova, Zwart and Boxma \cite{EG} as
\begin{equation}\label{1.2}
p(t)\sim C'\,e^{-A'\,t}\;\;\;(G=D),\;\;t\rightarrow\infty
\end{equation}
where the constants $A'$ and $C'$ are explicitly characterized in \cite{EG}. Comparing (\ref{1.1}) and (\ref{1.2}) we see that the tail behaviors of the $M/M/1$-PS and $M/D/1$-PS models are quite different.

In this paper we will study both the conditional sojourn time density $p(t|x)$ and the unconditional density $p(t)$ in the $M/G/1$-PS model. As in \cite{ZH} we shall consider various asymptotic limits, such as $x$ and/or $t\rightarrow\infty$ with a fixed $\rho<1$, and $1-\rho=\epsilon\rightarrow 0^+$ with space and time scaled by $\epsilon$. Here $\rho=\lambda\, m_1$ where $m_1=\int_0^\infty x\,b(x)dx$ is the mean service time. We consider service densities $b(x)$ that have ``thin tails", with decay that is at least exponential as $x\rightarrow \infty$. We shall show that the basic asymptotic structure of the conditional density $p(t|x)$ is essentially independent of the service density (though the formulas do depend on the Laplace transform of $b(x)$). In contrast, the unconditional density is highly dependent on the tail behavior of $b(x)$. We shall make specific assumptions on this tail, first assuming that
\begin{equation}\label{1.3}
b(x)\sim M\,x^q\,e^{-N\,x^r},\;\;r\geq 1,
\end{equation}
where $M,N$ ($>0$) and $q$ are constants. Thus (\ref{1.3}) allows for roughly exponential or even thinner tails, such as a Gaussian. Then we shall discuss ``zero-tail" service densities, by assuming that $b(x)$ has support for $0\leq x\leq A$ (e.g., $b(x)=1/A$ corresponds to uniformly distributed service times). In the zero tail case the behavior of $p(t)$ as $t\rightarrow \infty$ and $\rho<1$ is determined by the behavior of $b(x)$ near the upper limit of its support, and we will assume that
\begin{equation}\label{1.4}
b(x)\sim \alpha_*\,(A-x)^{\nu_*-1},\;\;x\uparrow A,
\end{equation}
where $\alpha_*$ and $\nu_*$ are positive constants.

We will obtain a wide variety of tail behaviors of $p(t)$ as $t\rightarrow\infty$ for the general $M/G/1$-PS model, that are different from either (\ref{1.1}) or (\ref{1.2}). We shall also identify the class of service densities that lead to purely exponential tails, such as $G=D$ in (\ref{1.2}).

We mention some related work on various PS models. Ramaswami \cite{RA} studied the $G/M/1$-PS queue and obtained explicit results for the unconditional moments of the sojourn time. Various asymptotic properties of the conditional and unconditional moments and distribution were derived in \cite{KN}. The $G/G/1$-PS model has not been analyzed exactly, but some approximations are discussed in Sengupta \cite{SEapp} and the tail exponent of the unconditional sojourn time density was derived by Mandjes and Zwart \cite{MA}. Specifically, in \cite{MA} the authors characterized the limit $t^{-1}\log[p(t)]\rightarrow -A_0$ as $t\rightarrow\infty$, assuming that the arrival and service densities have at least exponential tails. In \cite{ZW} Zwart and Boxma analyze the $M/G/1$-PS queue with heavy tails, where the service density has algebraic or sub-exponential behavior as $x\rightarrow\infty$ (thus $N=0, r=0, q<-1$ in (\ref{1.3}), or $0<r<1$).

For PS models one is also interested in the sojourn time conditioned on the number of other customers in the system when the tagged customer arrives. The conditional sojourn time for the $M/M/1$-PS model, conditioned on this number rather than the service time $x$, was studied by Sengupta and Jagerman \cite{SEcon} and Guillemin and Boyer \cite{GU}. A good recent survey of sojourn time asymptotics in PS queues is in Borst, N\'u\~nez-Queija and Zwart \cite{BONU}.

In this paper the main methods used are for the asymptotic expansion of integrals, such as the Laplace and saddle point methods, and good general references are the books of Bleistein and Handelsman \cite{BL} and of Wong \cite{WO}.

The remainder of the paper is organized as follows. In Section 2 we summarize and briefly discuss our main results (see Theorems 2.1--2.5). In Section 3 we derive the results for $p(t|x)$, for moderate traffic intensities $\rho<1$. In Section 4 we consider $p(t|x)$ for $\rho\uparrow 1$, and various scalings of space and time. We remove the condition on $x$ in Section 5, treating both $\rho<1$ and $\rho\approx 1$, and here we make the assumptions (\ref{1.3}) or (\ref{1.4}).

\section{Summary of results}
We assume that customers arrive according to a Poisson process with rate $\lambda$, at a single processor-sharing server. The customers' random service requests are i.i.d. random variables with density function $b(y)$, and Laplace-Stieltjes transform $\hat{b}(\tau)=\int_0^\infty e^{-\tau y}\,b(y)dy$. We assume that $\hat{b}(\tau)$ is an analytic function of $\tau$ for $\Re(\tau)>-\epsilon_0$ for some $\epsilon_0>0$. Thus all the moments of the service time are finite, and we set
\begin{equation}\label{moment}
m_k=\int_0^\infty y^k\,b(y)\,dy=(-1)^k\,\hat{b}^{(k)}(\tau)|_{\tau=0}\;\;\;(k\in\mathbb{N}).
\end{equation}
In particular we denote the service rate by $\mu$ where
\begin{equation*}\label{m1}
\frac{1}{\mu}=m_1=\int_0^\infty y\,b(y)\,dy=-\hat{b}^\prime(0).
\end{equation*}
Then the traffic intensity is defined as $\rho=\lambda/\mu$ and we assume that $\rho<1$.

The Laplace transform of the equilibrium sojourn time distribution, conditioned on the tagged customer (or job) requiring $x$ units of service, was derived by Ott \cite{OT} and Yashkov \cite{YAproc}, \cite{YAmath}, who obtained
\begin{equation}\label{Ee-sV}
\mathbf{E}[e^{-s \mathbf{V}(x)}]=\frac{1-\rho}{(1-\rho)\,G_1(s,x)+s\,G_3(s,x)},
\end{equation}
where
\begin{equation*}\label{G1}
\int_0^\infty{e^{-\tau x}\,G_1(s,x)}dx=\frac{\tau-\lambda\,(1-\hat{b}(\tau))}{\tau\big[\tau-s-\lambda\,(1-\hat{b}(\tau))\big]},
\end{equation*}
and
\begin{equation*}\label{G3}
\int_0^\infty{e^{-\tau x}\,G_3(s,x)}dx=\frac{\rho\,[\tau-\mu\,(1-\hat{b}(\tau))]}{\tau^2\,\big[\tau-s-\lambda\,(1-\hat{b}(\tau))\big]}.
\end{equation*}

Thus, the Laplace transform of the denominator in the right-hand side of (\ref{Ee-sV}) is
\begin{eqnarray}\label{f(tau,s)}
f(\tau;s) &=& \int_0^\infty\left[(1-\rho)\,G_1(s,x)+s\,G_3(s,x)\right]e^{-\tau x}dx\nonumber\\
          &=& \frac{(1-\rho)\,\tau^2-(1-\rho)\,\lambda\,(1-\hat{b}(\tau))\,\tau+s\,\rho\,\tau-s\,\lambda\,(1-\hat{b}(\tau))}{\tau^2\,\big[\tau-s-\lambda\,(1-\hat{b}(\tau))\big]}.
\end{eqnarray}
Taking the inverse Laplace transform of (\ref{f(tau,s)}), (\ref{Ee-sV}) becomes
\begin{equation}\label{Ee-sV2}
\mathbf{E}[e^{-s \mathbf{V}(x)}]=\frac{1-\rho}{\frac{1}{2\pi i}\int_{Br_\tau}e^{\tau x}f(\tau;s)\,d\tau},
\end{equation}
where $Br_{\tau}$ is a vertical contour in the complex $\tau$-plane, on which $\Re(\tau)>0$.

By taking the inverse Laplace transform of (\ref{Ee-sV2}), the
probability density of the sojourn time, conditioned on service
time $x$, is
\begin{eqnarray}\label{ptx_main}
p(t|x) &=& \frac{1}{2\pi i}\int_{Br_s}e^{st}\,\mathbf{E}[e^{-s \mathbf{V}(x)}]\,ds\nonumber\\
       &=& \frac{1-\rho}{2\pi i}\int_{Br_s}e^{st}\left[\frac{1}{2\pi i}\int_{Br_\tau}e^{\tau x}\,f(\tau;s)\,d\tau\right]^{-1}ds.
\end{eqnarray}

\noindent Here $Br_s$ is a vertical contour in the complex $s$-plane with $\Re(s)>0$. Note that $p(t|x)$ will in general have a probability mass along $t=x$. Analyzing the integral (\ref{ptx_main}), we obtain the following expansions for $p(t|x)$, valid on different space and time scales.

\begin{theorem} \label{th2.1}
For a fixed $\rho<1$, the conditional sojourn time density has the following asymptotic expansions:

\begin{enumerate}

\item $x\rightarrow \infty$, $t-x\rightarrow 0^+$ with $x(t-x)^\nu=O(1)$, assuming that
\begin{equation*}\label{by_in_th2.1}
b(y)\sim \alpha\, y^{\nu-1}, \textrm{  as  } y\rightarrow 0\;\;(\alpha, \nu>0),
\end{equation*}
\begin{eqnarray}\label{th2.1.1}
p(t|x) &\sim& \frac{1-\rho}{2\pi i}\int_{Br_s}e^{-\lambda\, x}\,e^{s\,(t-x)}\exp\Big[\frac{\lambda\,\alpha\,\Gamma(\nu)}{s^\nu}\,x\Big]ds\nonumber\\
       &=& (1-\rho)\,\delta(t-x)\,e^{-\lambda\, x}\nonumber\\
       && +(1-\rho)\,e^{-\lambda\, x}\sum_{m=1}^\infty\frac{[\lambda\,\alpha\,\Gamma(\nu)\,x]^m \,(t-x)^{\nu m-1}}{m!\;\Gamma(\nu m)}.
\end{eqnarray}

\item $x$,$t\rightarrow \infty$, $1<t/x<\infty$,
\begin{equation}\label{th2.1.2}
p(t|x)\sim\frac{(1-\rho)\,\tau_*^2\,\big[1+\lambda\,\hat{b}^\prime(\tau_*)\big]^{5/2}}{s_*^2\,\sqrt{2\pi x \lambda\,\hat{b}^{\prime\prime}(\tau_*)}}\,e^{s_*t}\,e^{-\tau_*x},
\end{equation}
where $s_*=s_*(t/x)$ and $\tau_*=\tau_*(t/x)$ satisfy the equations
\begin{equation}\label{th2.1.2*}
-\hat{b}^\prime(\tau_*)=\int_0^\infty e^{-\tau_*y}\,y\,b(y)\,dy=\frac{t-x}{\lambda\, t},\textrm{   } s_*=\tau_*-\lambda\,(1-\hat{b}(\tau_*)).
\end{equation}

\item $x,t\rightarrow\infty$, $t/x^2=O(1)$,
\begin{eqnarray*}\label{th2.1.3}
p(t|x)&\sim& \frac{(1-\rho)\,\tau_0^2\,\sqrt{\lambda\,\hat{b}^{\prime\prime}(\tau_0)}}{2^{3/2}\,\pi^{1/2}\,s_0^2\;t^{5/2}}\,e^{s_0\,t}\,e^{-\tau_0\,x}\sum_{n=0}^\infty\exp\Big\{-\frac{(2n+1)^2\,x^2}{\lambda\,\hat{b}^{\prime\prime}(\tau_0)\,t}\Big\}\nonumber\\
&&\times\bigg(\frac{2\,(2n+1)^2\,x^2}{\lambda\,\hat{b}^{\prime\prime}(\tau_0)}-2\,t\bigg),
\end{eqnarray*}
where $s_0=s_*(\infty)$ and $\tau_0=\tau_*(\infty)$ satisfy:
\begin{equation}\label{th2.1.3*}
\lambda\,\hat{b}^\prime(\tau_0)=-1,\textrm{    }s_0=\tau_0-\lambda\,(1-\hat{b}(\tau_0)).
\end{equation}

\item $x=O(1)$, $t\rightarrow\infty$,
\begin{equation}\label{th2.1.4}
p(t|x)\sim(1-\rho)\,J(x)\,e^{s_c(x)\,t},
\end{equation}
where
\begin{equation}\label{th2.1.4_J}
J(x)=\frac{d}{ds}\left(\frac{1}{2\pi i}\int_{Br_\tau}f(\tau;s)\,e^{\tau x}\,d\tau\right)\bigg\arrowvert_{s=s_c(x)},
\end{equation}
and $s_c(x)$($<0$) is the maximal real solution of
\begin{equation}\label{th2.1.4_Sc}
\frac{1}{2\pi i}\int_{Br_\tau}f(\tau;s_c)\,e^{\tau x}\,d\tau=0
\end{equation}
or
\begin{equation*}\label{th2.1.4_Sc(2)}
1-\rho+x\,s_c+\frac{s_c^2}{2\pi i}\int_{Br_\tau}\frac{e^{\tau x}}{\tau^2\,\big[\tau-s_c-\lambda\,(1-\hat{b}(\tau))\big]}\,d\tau=0.
\end{equation*}

\end{enumerate}

\end{theorem}

The result in Case 4 was also recently derived by Yashkov \cite{YAon}, who characterized $s_c(x)$ and $J(x)$ in a different form.

In the asymptotic matching region between Cases 3 and 4, we have
\begin{equation}\label{matching_th2.1}
p(t|x)\sim\frac{(1-\rho)\,\pi^2\,\lambda^2\,\tau_0^2\,(\hat{b}^{\prime\prime}(\tau_0))^2}{2\,s_0^2\,x^3}\,e^{-\tau_0\,x}\,\exp\left\{s_0\,t+B\,t/x^2+C\,t/x^3\right\},
\end{equation}
where
\begin{equation}\label{th2.1_BC}
B=-\frac{\pi^2\lambda\,\hat{b}^{\prime\prime}(\tau_0)}{2},\textrm{      }C=-\frac{\pi^2\lambda\,[6\,\hat{b}^{\prime\prime}(\tau_0)+\tau_0\,\hat{b}^{\prime\prime\prime}(\tau_0)]}{3\,\tau_0}.
\end{equation}

\noindent Then result (\ref{matching_th2.1}) holds for $x,t\rightarrow\infty$ with $t=O(x^3)$. It can be extended to larger ranges of $t$, e.g., to $t=O(x^4)$, by including an additional factor of the form $\exp(D\,t/x^4)$.

We note that in Case 2, if $t/x\sim 1/(1-\rho)$, by (\ref{th2.1.2*}) we have $s_*(1/(1-\rho))=0$, $\tau_*(1/(1-\rho))=0$ and
\begin{equation*}\label{}
s_*\,t-\tau_*\,x=-\frac{(1-\rho)^3}{2\,\lambda\, m_2\, x}\,\Big(t-\frac{x}{1-\rho}\Big)^2+O\Big(\big(t-\frac{x}{1-\rho}\big)^3\Big).
\end{equation*}
Then the formula (\ref{th2.1.2}) simplifies to the Gaussian
\begin{equation*}
p(t|x)\approx \frac{(1-\rho)^{3/2}}{\sqrt{2\,\pi\,\lambda\, m_2\, x}}\exp\bigg[-\frac{(1-\rho)^3}{2\,\lambda\, m_2 \,x}\,\Big(t-\frac{x}{1-\rho}\Big)^2\bigg],
\end{equation*}
which gives the spread about the well known mean value $\mathbf{E}\left[\mathbf{V}(x)\right]=x/(1-\rho)$.

As a special case of the $M/G/1$-PS model, we consider the $M/E_k/1$-PS model, in which the Erlang service time density function is given by
\begin{equation}\label{Erlang pdf}
b(y)=\frac{(k\,\mu)^k\,y^{k-1}\,e^{-k\,\mu\, y}}{(k-1)!}\;\;(k\in\mathbb{N},y\geq0)
\end{equation}
and thus
\begin{equation}\label{Erlang LST}
\hat{b}(\tau)=\left(\frac{k\,\mu}{k\,\mu+\tau}\right)^k.
\end{equation}

\noindent Then we obtain the following more explicit results.

\begin{corollary} \label{coro2.1}
For the $M/E_k/1$-PS model with traffic intensity $\rho<1$, the conditional sojourn time density has the following asymptotic expansions:

\begin{enumerate}

\item $x\rightarrow \infty$, $t-x\rightarrow 0^+$ with $x\,(t-x)^k=O(1)$,
\begin{eqnarray}\label{cor2.1.1}
p(t|x)&\sim& \frac{1-\rho}{2\pi i}\int_{Br_s}e^{-\lambda\, x}e^{s\,(t-x)}\exp\Big[\frac{\lambda\,(k\,\mu)^k}{s^k}\,x\Big]ds\nonumber\\
      &=&\; (1-\rho)\,\delta(t-x)\,e^{-\lambda\, x}+\frac{(1-\rho)\,\lambda\,(k\,\mu)^k\, x\,(t-x)^{k-1}\,e^{-\lambda\, x}}{(k-1)!}\nonumber\\
       &\times& _0F_k\big([\textrm{ }];[1+\frac{1}{k},1+\frac{2}{k},...,2-\frac{1}{k},2];\lambda\,\mu^k\,x\,(t-x)^k\big).
\end{eqnarray}
Here $_0F_k([\textrm{ }];[b_1,b_2,...,b_k];z)$ is the generalized hypergeometric function.

\item $x$,$\,t\rightarrow \infty$, $1<t/x<\infty$,
\begin{eqnarray}\label{cor2.1.2}
p(t|x) &\sim& \frac{(1-\rho)\,k\,\tilde{\tau}_*^2\,[t-(\mu\, k)^{k+1}\,(t-x)]\sqrt{\mu\,(\lambda\,\mu^k)^{\frac{1}{k+1}}\left(\frac{t}{t-x}\right)^{\frac{k+2}{k+1}}\left(\frac{x}{t}\right)^3}}{t\,\tilde{s}_*^2\,\sqrt{2\,\pi\,(k+1)\,x}}\nonumber\\
       && \times e^{\tilde{s}_*\,t}\,e^{-\tilde{\tau}_*\,x},
\end{eqnarray}
where
\begin{equation}\label{cor2.1.2s*}
\tilde{s}_*=(\lambda\,\mu^k)^{\frac{1}{k+1}}\,\Big(\frac{1}{1-x/t}\Big)^{\frac{1}{k+1}}(k+1-x/t)-\lambda-\mu\, k,
\end{equation}

\begin{equation}\label{cor2.1.2tau*}
\tilde{\tau}_*=k\,\mu\bigg[\Big(\frac{\rho}{1-x/t}\Big)^{\frac{1}{k+1}}-1\bigg].
\end{equation}

\item $x,t\rightarrow\infty$, $t/x^2=O(1)$,
\begin{eqnarray}\label{cor2.1.3}
p(t|x)&\sim& \frac{(1-\rho)\,k\,\sqrt{k\,(k+1)}\,(1-\rho^{\frac{1}{k+1}})^2}{2\sqrt{2\,\pi\,\mu}\,\rho^{\frac{1}{2(k+1)}}\,\big[\rho-(k+1)\,\rho^{\frac{1}{k+1}}+k\big]^2\,t^{5/2}}\,e^{\tilde{s}_0 t}\,e^{-\tilde{\tau}_0 t}\nonumber\\
&&\times \sum_{n=0}^\infty\exp\Big\{-\frac{(2n+1)^2\,k\,(\lambda\,\mu^k)^{\frac{1}{k+1}}\,x^2}{2\,(k+1)\,t}\Big\}\nonumber\\
&&\times\Big[\frac{2\,k}{k+1}(2\,n+1)^2\,(\lambda\,\mu^k)^{\frac{1}{k+1}}\,x^2-2\,t\Big],
\end{eqnarray}
where
\begin{equation}\label{cor2.1.3s0}
\tilde{s}_0=(k+1)\,(\lambda\,\mu^k)^{\frac{1}{k+1}}-(k\,\mu+\lambda),
\end{equation}
\begin{equation}\label{cor2.1.3tau0}
\tilde{\tau}_0=k\,(\lambda\,\mu^k)^{\frac{1}{k+1}}-k\,\mu.
\end{equation}

\item $x=O(1),t\rightarrow\infty$,
\begin{equation}\label{cor2.1.4}
p(t|x)\sim\frac{(1-\rho)\,e^{\tilde{s}_c(x)\,t}}{\frac{d}{ds}\Big[\sum_{i=1}^{k+1}e^{\tau_i(s)\,x}\,R_i(s)\Big]\Big\arrowvert_{s=\tilde{s}_c}},
\end{equation}
where $\tau_i=\tau_i(s)\;(i=1,...,k+1)$ are the $k+1$ poles of $f(\tau;s)$ in (\ref{f(tau,s)}), with residues $R_i(s)=Res(f,\tau=\tau_i(s))$, and $\tilde{s}_c(x)$ is the maximal root of $\sum_{i=1}^{k+1}e^{\tau_i(s)\,x}\,R_i(s)=0$.

\end{enumerate}
\end{corollary}

In the matching region between Cases 3 and 4, we have
\begin{equation}\label{matching_cor2.1}
 \begin{aligned}
p(t|x)\sim&\; \frac{(1-\rho)\,(k+1)^2\,\pi^2\,(1-\rho^{\frac{1}{k+1}})^2}{2\,\mu^2\,\rho^{\frac{2}{k+1}}\,\big[\rho-(k+1)\,\rho^{\frac{1}{k+1}}+k\big]^2\,x^3}\\
&\times\exp\left\{k\,\mu\,(1-\rho^{\frac{1}{k+1}})\,x+s_0 \,t+\tilde{B}\,t/{x^2}+\tilde{C}\,t/{x^3}\right\},
 \end{aligned}
\end{equation}
where
\begin{equation}\label{cor2.1_B}
\tilde{B}=-\frac{\pi^2\,(k+1)}{2\,k\,\mu\,\rho^{\frac{1}{k+1}}},
\end{equation}
\begin{equation}\label{cor2.1_C}
\tilde{C}=\frac{\pi^2\,(k+1)\big[(k-4)\,\rho^{\frac{1}{k+1}}-(k+2)\big]}{3\,k^2\,\mu^2\,\rho^{\frac{2}{k+1}}\,\big(\rho^{\frac{1}{k+1}}-1\big)}.
\end{equation}

We note that for the $M/E_k/1$-PS model, $\nu=k$, $\alpha=(k\,\mu)^k/{(k-1)!}$, and $\tau_*$, $s_*$, $s_0$ and $\tau_0$ in Theorem 2.1 are explicitly computable, as given by $\tilde{\tau}_*$, $\tilde{s}_*$, $\tilde{s}_0$ and $\tilde{\tau}_0$ in Corollary 2.1.

We next consider the heavy traffic case, where $\lambda\uparrow\mu$. Letting $\epsilon\equiv1-\rho$ (thus $\epsilon\rightarrow 0^+$), we have the following results for general service time distributions.

\begin{theorem}\label{th2.2}
For $\rho=1-\epsilon$, where $\epsilon\rightarrow 0^+$, we let $t=T/\epsilon$ and $x=X/\epsilon$. The conditional sojourn time density of $M/G/1$-PS model has the following asymptotic expansions:

\begin{enumerate}

\item $x=O(1), t=O(1)$,
\begin{equation}\label{th2.2.1}
p(t|x)\sim\frac{\epsilon}{2\pi i}\int_{Br_s}\frac{e^{s\,t}}{s}\bigg[\frac{1}{2\pi i}\int_{Br_\tau}\frac{e^{\tau\, x}[\tau-\mu+\mu\,\hat{b}(\tau)]}{\tau^2\,[\tau-\mu+\mu\,\hat{b}(\tau)-s]}\,d\tau\bigg]^{-1}ds.
\end{equation}

\item $x=O(1),T=O(1)$,
\begin{equation}\label{th2.2.2}
 \begin{aligned}
p(t|x) \sim\; & \frac{\epsilon}{x}e^{-T/x}-\epsilon^2\bigg[\frac{\delta(T)}{x^2}+\frac{(T-2x)\,e^{-T/x}}{x^4}\bigg]\\
        & \times\bigg[Q_*(x)+\frac{1}{2\pi i}\int_\mathcal{C_-}\frac{e^{\tau\, x}}{\tau^2[\tau-\mu\,(1-\hat{b}(\tau))]}\,d\tau\bigg],
 \end{aligned}
\end{equation}
where
\begin{eqnarray}\label{th2.2.2_Q*}
Q_*(x)&=&\frac{1}{270\,\mu\, m^4_2}\Big[90\,m^3_2\,x^3+90\,m^2_2\,m_3\,x^2+15\,m_2\,(4\,m^2_3-3\,m_2\,m_4)\,x\nonumber\\
&&+(9\,m^2_2\,m_5+20\,m_3^3-30\,m_2\,m_3\,m_4)\Big]
\end{eqnarray}
and $m_i$ is the $i^{th}$ moment of the service time distribution, given by (\ref{moment}). The contour $\mathcal{C_-}$ can be taken as the imaginary axis in the $\tau$-plane, indented to the left of $\tau=0$, where the integrand has a pole of order 4.

\item $X,T=O(1),T-X\rightarrow 0^+$ with $T-X=T_*\,\epsilon^{1+1/\nu}=O(\epsilon^{1+1/\nu})$,  assuming that
\begin{equation*}\label{by_in_th2.2}
b(y)\sim \alpha\, y^{\nu-1}, \textrm{  as    } y\rightarrow 0 \textrm{   }(\alpha, \nu>0),
\end{equation*}
\begin{eqnarray}\label{th2.2.3}
p(t|x) &\sim& \frac{e^{-\mu\, X/\epsilon}\epsilon^{1-1/\nu}}{2\pi i}\int_{Br_S}e^{S\,T_*}\exp\Big[\frac{\mu\,\alpha\,\Gamma(\nu)}{S^\nu}\,X\Big]\,dS\nonumber\\
       &=& e^{-\mu \,X/\epsilon}\epsilon^{1-1/\nu}\Big[\delta(T_*)+\sum_{m=1}^\infty\frac{T_*^{\nu \,m-1}\,[\mu\,\alpha\,\Gamma(\nu)\,X]^m}{m!\;\Gamma(\nu\, m)}\Big].
\end{eqnarray}

\item $X=O(1),T=O(1)$, and $1<T/X<\infty$,
\begin{equation}\label{th2.2.4}
p(t|x)\sim\frac{\epsilon^{3/2}\,\hat{s}_*^2}{\sqrt{2\,\pi\,\mu \,T\,\hat{b}^{\prime\prime}(\hat{\tau}_*)}}\,\exp\bigg[\frac{T\,\hat{s}_*-X\,\hat{\tau}_*}{\epsilon}+T(\hat{\tau}_*-\hat{s}_*)\bigg],
\end{equation}
where $\hat{s}_*=\hat{s}_*(T/X)$ and $\hat{\tau}_*=\hat{\tau}_*(T/X)$ satisfy
\begin{equation}\label{th2.2.4*}
1+\mu\,\hat{b}^\prime(\hat{\tau}_*)=X/T, \textrm{    }\hat{s}_*=\hat{\tau}_*-\mu\,(1-\hat{b}(\hat{\tau}_*)).
\end{equation}

\item $X=\sqrt{\epsilon}\,Z=O(\sqrt{\epsilon}), T=O(1)$,
\begin{equation}\label{th2.2.5}
p(t|x)\sim \frac{2\sqrt{2}\,\epsilon^{3/2}}{\sqrt{\mu m_2\,\pi \,T}}\sum_{n=0}^\infty \exp\bigg[-\frac{(2n+1)^2\,Z^2}{2\,\mu m_2\, T}\bigg].
\end{equation}

\item $X=O(1),T=\Theta/\epsilon=O(\epsilon^{-1})$, we give the expansion in three different forms:
\begin{enumerate}
\item
\begin{equation}\label{th2.2.6a}
 \begin{aligned}
p(t|x) &\sim \frac{\epsilon^2}{\mu m_2\,\pi\, i}\,\exp\Big(\frac{X}{\mu m_2}-\frac{\Theta}{2\,\mu m_2}\Big)\\
       &\times\int_{Br_\xi}\frac{\sqrt{\xi}\,\exp(\frac{\Theta}{2\,\mu m_2}\,\xi-\frac{X}{\mu m_2}\sqrt{\xi})}{(1+\sqrt{\xi})^2-(1-\sqrt{\xi})^2\exp(-\frac{2\,X}{\mu m_2}\sqrt{\xi})}\,d\xi.
 \end{aligned}
\end{equation}

\item
\begin{equation}\label{th2.2.6b}
 \begin{aligned}
p(t&|x)\sim \frac{\epsilon^2}{\sqrt{\pi}} \sum_{n=0}^\infty \exp\bigg[\frac{2\,(n+1)\,X}{\mu m_2}-\frac{z_n^2}{4}\bigg]\\
&\times\sum_{l=0}^{2n}(-1)^l\,\frac{(2n)!}{l!\,(2n-l)!}\,2^{l+3/2}\,(\mu m_2)^{-\frac{l+3}{2}}\,\Theta^{\frac{l+1}{2}}\\
&\times\bigg[\frac{\mu m_2}{\Theta}\,D_{-l}(z_n)-2\sqrt{\frac{\mu m_2}{\Theta}}\,D_{-l-1}(z_n)+D_{-l-2}(z_n)\bigg].\\
 \end{aligned}
\end{equation}
Here $D_\nu(\cdot)$ is the parabolic cylinder function and $z_n=\frac{(2n+1)X+\Theta}{\sqrt{\mu m_2\,\Theta}}$.

\item
\begin{equation}\label{th2.2.6c}
p(t|x)\sim \epsilon^2\sum_{n=1}^\infty e^{s_d(v_n)\,\Theta}G(v_n).
\end{equation}
Here $v_n=v_n(X)$ are the real positive roots of the equation

\begin{equation}\label{th2.2.6c_v}
\exp\Big[-\frac{2X}{\mu m_2}\,iv\Big]=\Big(\frac{1+iv}{1-iv}\Big)^2,
\end{equation}

\begin{equation}\label{th2.2.6c_G}
G(v_n)=\frac{2\,v_n^2\,\exp\big(\frac{X}{\mu m_2}\big)}{(v_n^2-1)X\cos\big(\frac{v_n\,X}{\mu m_2}\big)+\big[2\,v_n\, X+(v_n^2+1)\,\mu m_2\big]\sin\big(\frac{v_n\,X}{\mu m_2}\big)},
\end{equation}

\begin{equation}\label{th2.2.6c_Sd}
s_d(v_n)=B_1(v_n)+C_1(v_n)\,\epsilon+O(\epsilon^2),
\end{equation}
where
\begin{equation}\label{th2.2.6c_Bvn}
B_1(v_n)=-\frac{1+v_n^2}{2\,\mu m_2},
\end{equation}

\begin{eqnarray}\label{th2.2.6c_Cvn}
C_1(v_n)&=&\frac{(v_n^2+1)}{6\,\mu^2\,m_2^3\big[(v_n^2+1)\,X+2\mu m_2\big]}\Big[2\,\mu m_2\,(m_3-3\,\mu m_2^2)\nonumber\\
&&+\;\big(3\,\mu m_2^2\,v_n^2-3\,m_3\,v_n^2-3\,\mu m_2^2+m_3\big)\,X\Big].
\end{eqnarray}

\end{enumerate}
\end{enumerate}
\end{theorem}

For very large times, corresponding to $\Theta\gg 1$ (thus $t\gg\epsilon^{-2}$), we have
\begin{equation}\label{after th2.2}
p(t|x)\sim \epsilon^2 \,e^{s_d(v_1)\,\Theta}\,G(v_1),
\end{equation}
where $v_1=v_1(X)$ is the unique root of (\ref{th2.2.6c_v}) in the interval $(0,\mu m_2\,\pi/X)$. Then the first term in the sum in (\ref{th2.2.6c}) dominates.

For the $M/E_k/1$-PS model, we again get more explicit expressions.
\begin{corollary}
For the $M/E_k/1$-PS model in heavy traffic, we have the following expansions of the conditional sojourn time density.

\begin{enumerate}

\item $x=O(1),\,t=O(1)$,
\begin{equation}\label{cor2.2.1}
p(t|x)\sim\frac{\epsilon}{2\pi i}\int_{Br_s}\frac{e^{s\,t}}{s}\bigg[\frac{1}{2\pi i}\int_{Br_\tau}\frac{e^{\tau x}[(\tau-\mu)(k\mu+\tau)^k+\mu(k\mu)^k]}{\tau^2[(\tau-s-\mu)(k\mu+\tau)^k+\mu(k\mu)^k]}d\tau\bigg]^{-1}ds.
\end{equation}

\item $x=O(1),\,T=O(1)$,
\begin{equation}\label{cor2.2.2}
 \begin{aligned}
p(t|x) \sim\; &\frac{\epsilon}{x}e^{-T/x}-\epsilon^2\,\frac{x^2\,\delta(T)+(T-2x)\,e^{-T/x}}{x^4}\\
       & \times\bigg(\tilde{Q}_*(x)+\sum_{j=3}^{k+1}\frac{Q_j+k\,\mu}{(k+1)\,Q_j^3}\,e^{Q_j\,x}\bigg),
 \end{aligned}
\end{equation}
where
\begin{equation*}\label{cor2.2.2_Q*}
 \begin{aligned}
\tilde{Q}_*(x)=&\;\frac{1}{270(k+1)\,k^2\,\mu^2}\Big[90(k\,\mu)^3x^3+90(k+2)\,(k\,\mu)^2x^2\\
&+\;15k\,\mu\,(k^2+k-2)\,x-(k^3+9k^2+6k-16)\Big],
 \end{aligned}
\end{equation*}
and $Q_j$, $j=3,...,k+1$ are nonzero roots of
$$ (Q-\mu)(Q+k\,\mu)^k+\mu\,(k\,\mu)^k=0. $$
Note that $Q=0$ is a double root and the other roots have $\Re(Q)<0$. For example, if $k=2$, $Q_3=-3\,\mu$; if $k=3$, $Q_3=(-4+i\sqrt{2})\,\mu$ and $Q_4=(-4-i\sqrt{2})\,\mu$.

\item $X,T=O(1),\,T-X\rightarrow 0^+$ with $T-X=T_*\,\epsilon^{1+1/k}=O(\epsilon^{1+1/k})$,
\begin{eqnarray*}\label{cor2.2.3}
p(t|x)&\sim& \frac{e^{-\mu\, X/\epsilon}\epsilon^{1-1/\nu}}{2\pi i}\int_{Br_S}e^{S\,T_*}\exp\Big[\frac{\mu\,(k\,\mu)^k\,X}{S^k}\Big]\,dS\nonumber\\
&=&  e^{-\mu\, X/\epsilon}\,\epsilon^{1-1/k}\,\delta(T_*)+e^{-\mu\, X/\epsilon}\frac{\mu\,(k\,\mu)^k X\,T_*^{k-1}}{(k-1)!}\,\epsilon^{1-1/k}\nonumber\\
&\times& _0F_k\big([\textrm{ }];[1+\frac{1}{k},1+\frac{2}{k},...,2-\frac{1}{k},2];\mu^{k+1}X\,T_*^k\big).
\end{eqnarray*}

\item $X=O(1),\,T=O(1)$ and $1<T/X<\infty$,
\begin{eqnarray*}\label{cor2.2.4}
p(t|x) &\sim& \frac{\epsilon^{3/2}\,k^{5/2}\sqrt{\mu}\,\bigg[\left(\frac{T}{T-X}\right)^{\frac{1}{k+1}}-1\bigg]^2\left(\frac{X}{T-X}\right)^{\frac{k+2}{2(k+1)}}}{\sqrt{2\pi(k+1)X}\left(\frac{T}{X}\right)^{\frac{4k+3}{2(k+1)}}\left[\left(\frac{T-X}{T}\right)^{\frac{1}{k+1}}\left(\frac{kT}{T-X}+1\right)-(k+1)\right]^2}\nonumber\\
        && \times \exp\bigg\{\frac{X}{\epsilon}\,\mu\,\Big[(k+1)\,\frac{T}{X}\,\Big(\Big(\frac{T-X}{X}\Big)^{\frac{k}{k+1}}-1\Big)+k\Big] \bigg\}\nonumber\\
        && \times \exp\bigg\{X\,\mu\,\Big[\frac{T}{X}-\Big(\frac{T-X}{X}\Big)\Big(\frac{T}{T-X}\Big)^{\frac{1}{k+1}}  \Big]\bigg\}.
\end{eqnarray*}

\item $X=\sqrt{\epsilon}\,Z=O(\sqrt{\epsilon}),\,T=O(1)$,
\begin{equation*}\label{cor2.2.5}
p(t|x)\sim \epsilon^{3/2}\sqrt{\frac{8\,k\,\mu}{(k+1)\,\pi\, T}}\sum_{n=0}^\infty \exp\bigg[-\frac{(2n+1)^2\,k\,\mu\, Z^2}{2\,(k+1)\,T}\bigg].
\end{equation*}

\item $X=O(1),\,T=\Theta/\epsilon=O(\epsilon^{-1})$, we have the following three different forms of the expansion:
\begin{enumerate}
\item
\begin{equation*}\label{cor2.2.6a}
 \begin{aligned}
p(t|&x) \sim \frac{\epsilon^2\, \mu\, k}{(k+1)\,\pi\, i}\exp\bigg[\frac{\mu\, k\, X}{k+1}-\frac{\mu\, k\,\Theta}{2\,(k+1)}\bigg]\\
&\times\int_{Br_\xi}\frac{\sqrt{\xi}\,\exp\left(\frac{\mu\, k\,\Theta}{2(k+1)}\xi-\frac{\mu\, k\, X}{k+1}\sqrt{\xi}\right)}{(1+\sqrt{\xi})^2-(1-\sqrt{\xi})^2\exp\Big(-\frac{2\,\mu\, k\, X}{k+1}\sqrt{\xi}\Big)}\,d\xi.\\
 \end{aligned}
\end{equation*}

\item
\begin{eqnarray*}\label{cor2.2.6b}
p(t|x)&\sim& \frac{\epsilon^2}{\sqrt{\pi}} \sum_{n=0}^\infty \exp\bigg[\frac{2\,k\,\mu\,(n+1)\,X}{k+1}-\frac{\tilde{z}_n^2}{4}\bigg]\sum_{l=0}^{2n}(-1)^l\frac{(2n)!}{l!\,(2n-l)!}\nonumber\\
&&\times\;2^{l+3/2}\Big(\frac{k}{k+1}\Big)^{\frac{l+3}{2}}\,\mu^{l/2-1}\,\Theta^{\frac{l+1}{2}}\bigg[\frac{(k+1)\,\mu}{k\,\Theta}\,D_{-l}(\tilde{z}_n)\nonumber\\
&& -2\sqrt{\frac{(k+1)\,\mu}{k\,\Theta}}\,D_{-l-1}(\tilde{z}_n)+\mu^2\,D_{-l-2}(\tilde{z}_n)\bigg],
\end{eqnarray*}
where $\tilde{z}_n=\sqrt{\frac{k\,\mu}{(k+1)\,\Theta}}\big[\Theta+(2n+1)\,X\big]$.

\item
\begin{equation*}\label{cor2.2.6c}
p(t|x)\sim\epsilon^2\sum_{n=1}^\infty e^{s_d(\tilde{v}_n)\,\Theta}\,\tilde{G}(\tilde{v}_n).
\end{equation*}
Here $\tilde{v}_n=\tilde{v}_n(X)$ are the real positive roots of the equation
\begin{equation}\label{cor2.2.6c_v}
\exp\bigg[-\frac{2\,k\,\mu\, X}{k+1}i\tilde{v}\bigg]=\left(\frac{1+i\tilde{v}}{1-i\tilde{v}}\right)^2,
\end{equation}
\begin{equation*}\label{cor2.2.6c_G}
\tilde{G}(\tilde{v}_n)=\frac{2\,\tilde{v}_n^2\,\exp\left(\frac{\mu\, k\,X}{k+1}\right)}{(\tilde{v}_n^2-1)\,X\,\cos\left(\frac{k\,\mu\,\tilde{v}_n\,X}{k+1}\right)+\left[2\,\tilde{v}_n\, X+\frac{(k+1)\,(\tilde{v}_n^2+1)}{k\,\mu\,\tilde{v}_n}\right]\sin\left(\frac{k\,\mu\,\tilde{v}_n\,X}{k+1}\right)},
\end{equation*}
\begin{equation*}\label{cor2.2.6c_Sd}
s_d(\tilde{v}_n)=\tilde{B}_1(\tilde{v}_n)+\tilde{C}_1(\tilde{v}_n)\,\epsilon+O(\epsilon^2),
\end{equation*}
where
\begin{equation*}\label{cor2.2.6c_Btilde}
\tilde{B}_1(\tilde{v}_n)=-\frac{k\,\mu\,(\tilde{v}_n^2+1)}{2(k+1)},
\end{equation*}
\begin{equation*}\label{cor2.2.6c_Ctilde}
\tilde{C}_1(\tilde{v}_n)=-\frac{k\,\mu\,(\tilde{v}_n^2+1)}{6(k+1)^2}\bigg[(2k+1)-\frac{2\,(k-1)\,k\,\mu\, X\,\tilde{v}_n^2}{k\,\mu\,(\tilde{v}^2_n+1)\,X+2(k+1)}\bigg].
\end{equation*}

\end{enumerate}
\end{enumerate}
\end{corollary}

For very large times with $\Theta\gg 1$, we again have
\begin{equation*}\label{after cor2.2}
p(t|x)\sim \epsilon^2\, e^{s_d(\tilde{v}_1)\,\Theta}\,\tilde{G}(\tilde{v}_1),
\end{equation*}
where $\tilde{v}_1=\tilde{v}_1(X)$ is the unique root of (\ref{cor2.2.6c_v}) in the interval $(0,\frac{(k+1)\,\pi}{k\,X})$.

The asymptotics of the conditional sojourn time distribution are generally dependent on the service density only through equations such as (\ref{th2.1.2*}) or (\ref{th2.1.3*}), which involve the Laplace transform. This is true for all the scales in Theorem 2.1 except Case 1, where $t\approx x$ and the expansion depends on the local behavior of the service density $b(y)$ as $y\rightarrow 0^+$. We now examine the unconditional sojourn time density, defined by
\begin{equation}\label{p(t)}
p(t)=\int_0^t b(x)\,p(t|x)dx.
\end{equation}
The structure of $p(t)$ is highly dependent on the behavior of $b(x)$ for $x\rightarrow\infty$, and we need to make specific assumptions on the tail of the service distribution. We first assume that the service time density function behaves as
\begin{equation}\label{b(y)_exp}
b(y)\sim M\,y^q\,e^{-N\,y^r},\;y\rightarrow\infty\;\;(M,N>0, \,1\leq r\leq 2).
\end{equation}

We remove the condition on $x$ by using the results in (\ref{matching_th2.1}) (if $\rho<1$) or (\ref{after th2.2}) (if $\rho\sim 1$), and evaluate asymptotically the integral in (\ref{p(t)}). We thus derive the following results for the unconditional sojourn time density as $t\rightarrow \infty$, from which we realize the variety of possible tail behaviors for the $M/G/1$-PS model.

\begin{theorem}
As $t\rightarrow\infty$, the unconditional sojourn time density has the following asymptotic expansions, assuming the tail behavior in (\ref{b(y)_exp}), and that $\rho$ is fixed with $\rho<1$.

\begin{enumerate}
\item If $r=1$,
\begin{equation}\label{th2.3.1}
p(t)\sim\alpha_1\,t^{\frac{2q-5}{6}}\exp(s_0\,t-\gamma_1\,t^{1/3}),
\end{equation}
where
\begin{eqnarray}\label{alpha1}
\alpha_1 &=& \frac{M(1-\rho)\tau_0^2}{\sqrt{6}\,s_0^2\,(N+\tau_0)^{\frac{q-1}{3}}}\,\pi^{\frac{4q+5}{6}}\,\big[\lambda\,\hat{b}^{\prime\prime}(\tau_0)\big]^{\frac{2q+7}{6}}\nonumber\\
         && \times\exp\bigg\{-\frac{(N+\tau_0)\big[6\,\hat{b}^{\prime\prime}(\tau_0)+\tau_0\,\hat{b}^{\prime\prime\prime}(\tau_0)\big]}{3\,\tau_0\,\hat{b}^{\prime\prime}(\tau_0)}\bigg\},\nonumber
\end{eqnarray}
\begin{equation*}\label{gamma1}
\gamma_1=\frac{3}{2}(N+\tau_0)^{2/3}\,\big[\pi^2\,\lambda\,\hat{b}^{\prime\prime}(\tau_0)\big]^{1/3}
\end{equation*}
and $s_0$ and $\tau_0$ are computed from (\ref{th2.1.3*}).

\item If $1<r<2$,
\begin{equation}\label{th2.3.2}
p(t)\sim\alpha_2\,t^{\frac{2q-r-4}{2(r+2)}}\exp(s_0\,t-\phi(\xi_*,t)).
\end{equation}
Here
\begin{equation}\label{th2.3.2_phi}
\phi(\xi,t)=\Big(N\xi^r-\frac{B}{\xi^2}\Big)\,t^{\frac{r}{r+2}}+\tau_0\,\xi\, t^{\frac{1}{r+2}}-\frac{C}{\xi^3}\,t^{\frac{r-1}{r+2}},
\end{equation}
\begin{equation*}\label{alpha2}
\alpha_2=\frac{(1-\rho)\,M\,\tau_0^2}{\sqrt{2(r+2)}\,s_0^2}(N\,r)^{\frac{1-q}{r+2}}\,\pi^{\frac{3}{2}+\frac{2(q-1)}{r+2}}\big[\lambda\,\hat{b}^{\prime\prime}(\tau_0)\big]^{\frac{3}{2}+\frac{q-1}{r+2}},
\end{equation*}
where $B$ and $C$ are given by (\ref{th2.1_BC}), and $\xi_*=\xi_*(t)$ satisfies
\begin{equation}\label{th2.3.2_phi'}
\phi_\xi(\xi,t)=\Big(N\,r\,\xi^{r-1}+\frac{2B}{\xi^3}\Big)\,t^{\frac{r}{r+2}}+\tau_0\, t^{\frac{1}{r+2}}+\frac{3C}{\xi^4}\,t^{\frac{r-1}{r+2}}=0.
\end{equation}
We have the following expansion(s) of $\xi_*$  and $\phi(\xi_*,t)$ as $t\rightarrow\infty$:
\begin{eqnarray*}\label{xi*}
\xi_* &\sim& \bigg[\frac{\pi^2\,\lambda\,\hat{b}^{\prime\prime}(\tau_0)}{Nr}\bigg]^{\frac{1}{r+2}}+\frac{\tau_0\,[\pi^2\,\lambda\,\hat{b}^{\prime\prime}(\tau_0)]^{\frac{2-r}{r+2}}}{(r+2)(N\,r)^{\frac{4}{r+2}}}\,t^{\frac{1-r}{r+2}}\\
 &&+ \left\{
   \begin{aligned}
&\frac{(5-r)\,\tau_0^2\,[\pi^2\,\lambda\,\hat{b}^{\prime\prime}(\tau_0)]^{\frac{3-2r}{r+2}}}{2(r+2)^2(Nr)^{\frac{7}{r+2}}}\,t^{\frac{2(1-r)}{r+2}},  &\; 1<r\leq 3/2,  \\
   &\frac{6\,\hat{b}^{\prime\prime}(\tau_0)+\tau_0\,\hat{b}^{\prime\prime\prime}(\tau_0)}{(r+2)\,\tau_0\,\hat{b}^{\prime\prime}(\tau_0)}\,t^{-\frac{1}{r+2}},  &\; 3/2<r<2,  \\
  \end{aligned}
  \right.
\end{eqnarray*}

\begin{eqnarray*}\label{phi(xi*)}
\phi(\xi_*,t) &\sim& \frac{r+2}{2\,r}(Nr)^{\frac{2}{r+2}}\Big[\pi^2\,\lambda\,\hat{b}^{\prime\prime}(\tau_0)\Big]^{\frac{r}{r+2}}\,t^{\frac{r}{r+2}}+\tau_0\,\bigg[\frac{\pi^2\,\lambda\,\hat{b}^{\prime\prime}(\tau_0)}{Nr}\bigg]^\frac{1}{r+2}t^{\frac{1}{r+2}}\\
 &+& \left\{
   \begin{aligned}
&\frac{3\,\tau_0^2}{2(r+2)}(Nr)^{-\frac{4}{r+2}}\Big[\pi^2\,\lambda\,\hat{b}^{\prime\prime}(\tau_0)\Big]^{\frac{2-r}{r+2}}\,t^{\frac{2-r}{r+2}},  &\; 1<r\leq 3/2,  \\
   &\frac{\pi^2\,\lambda\,\big[6\,\hat{b}^{\prime\prime}(\tau_0)+\tau_0\,\hat{b}^{\prime\prime\prime}(\tau_0)\big]}{3\,\tau_0}\bigg[\frac{Nr}{\pi^2\,\lambda\,\hat{b}^{\prime\prime}(\tau_0)}\bigg]^{\frac{3}{r+2}}\,t^{\frac{r-1}{r+2}},  &\; 3/2<r<2.  \\
  \end{aligned}
  \right.
\end{eqnarray*}

\item If $r=2$,
\begin{equation}\label{th2.3.3}
p(t)\sim\alpha_3\,t^{\frac{q-3}{4}}\exp(s_0\,t-\gamma_3\, t^{1/2}-\delta_3\, t^{1/4}),
\end{equation}
where
\begin{equation*}\label{alpha3}
\alpha_3=\frac{(1-\rho)\,M\,\pi^{q/2+1}\,\tau_0^2\,\big[\lambda\,\hat{b}^{\prime\prime}(\tau_0)\big]^{\frac{q+5}{4}}}{2^{\frac{q+5}{4}}\,s_0^2 \,N^{\frac{q-1}{4}}}\exp\bigg[\frac{(3NC-\tau_0\,B)^2}{16NB^2}-\frac{ND}{B}\bigg],
\end{equation*}
$$ \gamma_3=\pi\sqrt{2\,N\,\lambda\,\hat{b}^{\prime\prime}(\tau_0)},$$
$$ \delta_3=-\frac{N\,C+\tau_0\,B}{N^{1/4}\,|B|^{3/4}}, $$
and
\begin{eqnarray*}
D&=&\frac{\pi^2\,\lambda}{72\,\tau_0^2\,\hat{b}^{\prime\prime}(\tau_0)}\Big\{\big[3\,\pi^2\,\hat{b}^{\prime\prime}(\tau_0)\,\hat{b}^{(4)}(\tau_0)-(5\pi^2+12)(\hat{b}^{\prime\prime\prime}(\tau_0))^2\big]\,\tau_0^2\\
&& -144\,\hat{b}^{\prime\prime}(\tau_0)\,\hat{b}^{\prime\prime\prime}(\tau_0)\,\tau_0-432\,(\hat{b}^{\prime\prime}(\tau_0))^2\Big\}.
\end{eqnarray*}

\end{enumerate}
\end{theorem}

For the $M/E_k/1$-PS model, we note that the parameters in the tail of the service time density in (\ref{b(y)_exp}) are $M=(k\,\mu)^k/{(k-1)!}$, $N=k\,\mu$, $q=k-1$ and $r=1$. Thus Case 1 of Theorem 2.3 applies, with $\alpha_1$ and $\gamma_1$ given by
\begin{eqnarray*}\label{alpha1_Ek}
\alpha_1 &=& \frac{\pi^{\frac{4k+1}{6}}\,k^{\frac{2k+11}{6}}\,(k+1)^{\frac{2k+5}{6}}\,\mu^{\frac{2k-1}{6}}\,(1-\rho)\,(1-\rho^{\frac{1}{k+1}})^2}{\sqrt{6}\,(k-1)!\,\rho^{\frac{4k+1}{6k+6}}\,\big[\rho-(k+1)\,\rho^{\frac{1}{k+1}}+k\big]^2}\\
&&\times\exp\bigg[\frac{k+2-(k-4)\,\rho^{\frac{1}{k+1}}}{3\,(1-\rho^{\frac{1}{k+1}})}\bigg],
\end{eqnarray*}
\begin{equation*}\label{gamma1_Ek}
\gamma_1=\frac{3}{2}\,\pi^{2/3}\,[k\,(k+1)]^{1/3}\,(\mu^k\,\lambda)^{\frac{1}{3(k+1)}},
\end{equation*}
and $s_0$ is given by $\tilde{s}_0$ in (\ref{cor2.1.3s0}).

We next consider the unconditional distribution in the heavy
traffic limit. We find that the expansion of $p(t)$ in
(\ref{p(t)}) is different for the time scales $t=O(1)$,
$t=O(\epsilon^{-1})$ and $t=O(\epsilon^{-r-2})$. For $t=O(1)$ we
can use (\ref{cor2.2.1}) to evaluate the integral in (\ref{p(t)}),
but no further simplification is possible. Below we give the
results on the large times scales, and we note that most of the
probability mass occurs where $t=T/\epsilon=O(\epsilon^{-1})$.

\begin{theorem}
For $\rho=1-\epsilon$, where $\epsilon\rightarrow 0^+$, the unconditional sojourn time density has the following asymptotic expansions.
\begin{enumerate}
\item $t=T/\epsilon=O(\epsilon^{-1})$,
\begin{equation}\label{th2.4.1}
p(t)\sim\epsilon\int_0^\infty\frac{e^{-T/x}}{x}\,b(x)\,dx+\epsilon^2 \,S(T),
\end{equation}
where $S(T)$ is given by the double integral
\begin{equation*}\label{S(T)}
S(T)=-\int_0^\infty b(x)\,\frac{(T-2x)\,e^{-T/x}}{x^4}\bigg[\frac{1}{2\pi i}\int_{\mathcal{C_+}}\frac{e^{\tau \,x}}{\tau^2\,\big[\tau-\mu\,(1-\hat{b}(\tau))\big]}d\tau\bigg]dx.
\end{equation*}
Here the contour $\mathcal{C_+}$ is taken as the imaginary axis in the $\tau$-plane, indented to the right of $\tau=0$.

\item $t=\sigma_*/{\epsilon^{r+2}}=O(\epsilon^{-r-2})$ with $1\leq r<2$,
\begin{eqnarray}\label{th2.4.2}
p(t)&\sim& \frac{\sqrt{2\,\pi\,\mu m_2}\,M\,\hat{X}^q\,G(v_1(\hat{X}))\,\epsilon^{r/2-q}}{\sqrt{\sigma_*\,\Big[(v_1'(\hat{X}))^2+v_1(\hat{X})\,v_1''(\hat{X})\Big]+r\,(r-1)\,\mu m_2\,N\,X^{r-2}}}\nonumber\\
&&\times\exp\left[-\psi_0(\hat{X}(\sigma_*),\sigma_*)\,\epsilon^{-r}+C_1(v_1(\hat{X}))\,\sigma_*\,\epsilon^{1-r}\right].
\end{eqnarray}

Here $\hat{X}=\hat{X}(\sigma_*)$ satisfies
\begin{equation}\label{th2.4.2_sigma*}
\sigma_*=\frac{N\,r\,X^{r-1}\,\big[\mu m_2\,X\,(v_1^2(X)+1)+2\,\mu^2\,m_2^2\big]}{v_1^2(X)\,\big[v_1^2(X)+1\big]},
\end{equation}
$\psi_0(X,\sigma_*)$ is defined by
$$\psi_0(X,\sigma_*)=N\, X^r-B_1(v_1(X))\,\sigma_*,$$
$v_1=v_1(X)$ is the unique root of (\ref{th2.2.6c_v}) in the interval $(0,\mu m_2\,\pi/X)$, and $B_1(v_1)$ and $C_1(v_1)$ are given by (\ref{th2.2.6c_Bvn}) and (\ref{th2.2.6c_Cvn}).

\item $t=\sigma/{\epsilon^4}=O(\epsilon^{-4})$ with $r=2$,
\begin{eqnarray}\label{th2.4.3}
p(t)&\sim& \frac{\sqrt{2\,\pi\,\mu m_2}\,M\,\tilde{X}^q\,G(v_1(\tilde{X}))\,\epsilon^{1-q}}{\sqrt{\sigma\,    \big[(v_1'(\tilde{X}))^2+v_1(\tilde{X})\,v_1''(\tilde{X})\big]+2\mu m_2\,N}}\nonumber\\
&&\times\;\exp\bigg\{-\psi_0(\tilde{X}(\sigma),\sigma)\,\epsilon^{-2}+C_1(v_1(\tilde{X}))\,\sigma\,\epsilon^{-1}+D_1(v_1(\tilde{X}))\,\sigma\nonumber\\
&&+\;\frac{\mu m_2\,\big[C_1'(v_1(\tilde{X}))\,\sigma\big]^2}{4\,\mu m_2\,N+2\,\sigma\,\big[(v_1'(\tilde{X}))^2+v_1(\tilde{X})\,v_1''(\tilde{X})\big]}\bigg\}.
\end{eqnarray}
Here $\tilde{X}=\tilde{X}(\sigma)$ satisfies
\begin{equation*}\label{th2.4.3_sigma}
\sigma=\frac{2NX\big[\mu m_2\,X\,(v_1^2(X)+1)+2\,\mu^2\,m_2^2\big]}{v_1^2(X)\, 	\big[v_1^2(X)+1\big]},
\end{equation*}
and $D_1(v)$ is defined by
\begin{equation*}\label{th2.4.3_D1}
D_1(v)=\frac{(v^2+1)\,\big(d_0+d_1\,X+d_2\,X^2+d_3\,X^3\big)}{72\,\mu^3\, m_2^5\,\big[X\,(v^2+1)+2\,\mu m_2\big]^3},
\end{equation*}
where
\begin{equation*}\label{d0}
d_0=24\,\mu^3\,m_2^3\,\Big[-12\,\mu^2\,m_2^4+8\,\mu m_2^2\,m_3+(v^2+1)\,m_2\,m_4-(v^2+3)\,m_3^2\Big],
\end{equation*}
\begin{eqnarray*}\label{d1}
d_1&=&4\,\mu^2\,m_2^2\,\Big[36\,\mu^2\,(v^2-3)\,m_2^4-24\,\mu\,(4v^2-3)\,m_2^2\,m_3\\
&&+3\,(3\,v^4-2\,v^2+3)\,m_2\,m_4-(11\,v^4-42\,v^2+27)\,m_3^2\Big],
\end{eqnarray*}
\begin{eqnarray*}\label{d2}
d_2&=&2\,\mu m_2\Big[-36\,\mu^2\,(v^2+3)\,m_2^4-24\,\mu\,(v^4+2\,v^2-3)\,m_2^2\,m_3\\
&+&3\,(3\,v^6-7\,v^4-7\,v^2+3)\,m_2\,m_4-(13\,v^6-27\,v^4-45\,v^2+27)\,m_3^2\Big],
\end{eqnarray*}
\begin{eqnarray*}\label{d3}
d_3&=&(v^2+1)^2\,\Big[-36\,\mu^2\,m_2^4+24\,\mu m_2^2\,m_3-3\,(v^4-6\,v^2+1)\,m_2\,m_4\\
&&-(5\,v^4-18\,v^2+9)\,m_3^2\Big].
\end{eqnarray*}

\end{enumerate}
\end{theorem}

If the service density had even thinner tails, say with $r>2$ in (\ref{b(y)_exp}), we can easily extend Theorem 2.3. The main complication is that we would need further terms in the expansion of $s_c(x)$ as $x\rightarrow\infty$, which has the form $s_0+B/x^2+C/x^3+D/x^4+O(x^{-5})$. When $r=2$ the $O(x^{-4})$ term affects the leading term in the expansion of $p(t)$ as $t\rightarrow\infty$. The asymptotic evaluation of (\ref{p(t)}) involves balancing the factors $\exp(B\,t/x^2)$ and $\exp(-N\,x^r)$, which occurs when $x=O(t^{\frac{1}{r+2}})$. We can also extend Theorem 2.3 to more complicated tail behaviors of $b(y)$ of the form
\begin{equation*}\label{b(y)_expmore}
b(y)\sim M\,y^q\,\exp\big[-N\,y^r+N_1\,y^{r_1}+\cdots+N_l\,y^{r_l}\big]
\end{equation*}
where $0<r_l<r_{l-1}<\cdots<r_1<r$. This would be needed, for example, if $b(y)$ were a truncated Gaussian centered at some non-zero $y$ (then $r=2$ and $r_1=1$). We shall next consider densities with ``zero-tail", and these lead to different behaviors of $p(t)$.

Now we assume that the service time density function $b(y)$ has
finite support for $0\leq y\leq A$ and behaves as
\begin{equation}\label{b(y)_support}
b(y)\sim \alpha_*\,(A-y)^{\nu_*-1},\;\;y\uparrow A,\;\;(\alpha_*,\,\nu_*>0).
\end{equation}
The structure of the conditional sojourn time density $p(t|x)$ is the same as in Theorem 2.1, but the unconditional sojourn time density $p(t)$ is determined by the behavior of $b(y)$ near the upper limit $A$ of its support.

As $t\rightarrow\infty$ with fixed $\rho$ and $\rho<1$, we remove
the condition on $x$ by using the results in (\ref{th2.1.4}). For
the heavy traffic case with $\rho\sim 1$, we remove the condition
on $x$ by using (\ref{th2.2.2}) on the large time scale
$t=O(\epsilon^{-1})$. We thus have the following results for the
unconditional sojourn time density.

\begin{theorem}\label{th2.5}
The unconditional sojourn time density has the following asymptotic expansions, assuming the service density behavior in (\ref{b(y)_support}).
\begin{enumerate}
\item $t\rightarrow\infty$ with $\rho$ fixed and $\rho<1$:
\begin{equation}\label{th2.5.1}
p(t)\sim \frac{(1-\rho)\,\alpha_*\,\Gamma(\nu_*)}{[s_c'(A)\,t]^{\nu_*}}\,J(A)\,e^{s_c(A)\,t},
\end{equation}
where $J(x)$ is given by (\ref{th2.1.4_J}) and $s_c(x)$ is the maximal real solution of (\ref{th2.1.4_Sc}).

\item $t=T/\epsilon=O(\epsilon^{-1})$, $\rho=1-\epsilon$ with $\epsilon\rightarrow 0^+$:
\begin{equation}\label{th2.5.2}
p(t)= \epsilon\,\int_0^A b(x)\,\frac{e^{-T/x}}{x}\,dx+O(\epsilon^2).
\end{equation}

\end{enumerate}
\end{theorem}

For very large times with $T\rightarrow\infty$, (\ref{th2.5.2}) becomes\
\begin{equation}\label{after_th2.5}
p(t)\sim \epsilon\,\alpha_*\,\Gamma(\nu_*)\,A^{2\nu_*-1}\,T^{-\nu_*}\,e^{-T/A}.
\end{equation}
We can show that when (\ref{th2.5.1}) is expanded in the heavy traffic limit $\rho\uparrow 1$, we also obtain (\ref{after_th2.5}). We note that as $\rho\uparrow 1$, (\ref{th2.1.4_Sc}) shows that $s_c(A)\sim -\epsilon/A=O(\epsilon)$ and thus $s_c'(A)t \sim T/{A^2}$. The $O(\epsilon^2)$ term in (\ref{th2.5.2}) is the same as that in (\ref{th2.4.1}), except that the integral over $x$ is truncated at $x=A$.

If $b(A)$ is non-zero and finite then $\nu_*=1$ and (\ref{th2.5.1}) shows that $p(t)$ has an exponential tail with the additional algebraic factor $1/t$. This additional factor disappears only in the limit of $\nu_*\rightarrow 0$, but then (\ref{b(y)_support}) shows that $b(y)$ develops a probability mass at $y=A$. Our results show that $p(t)$ will have a purely exponential tail only if $b(y)$ consists of one (then $G=D$) or several point masses, or if $b(y)$ has a point mass at the maximum of its support ($y=A$) with all the remaining mass in the range $0\leq y\leq A_*$ with $A_*<A$.

\section{Brief derivations of the conditional sojourn time density for the case $\rho<1$}

In this section, we first give a brief derivation of the conditional sojourn time density for the $M/E_k/1$-PS model, with a fixed traffic intensity $\rho$ less than one. Then we sketch the derivations for the general service time density.

The Erlang service time density and its Laplace transform are given by (\ref{Erlang pdf}) and (\ref{Erlang LST}). Then (\ref{f(tau,s)}) has the following explicit form:
\begin{equation}\label{S3_f_Er}
   \begin{aligned}
&f(\tau;s)\\
&=\frac{\tau(1-\rho)\big[(\tau-\lambda)(\tau+k\,\mu)^k+\lambda(k\,\mu)^k\big]+s\rho\big[(\tau-\mu)(\tau+k\,\mu)^k+\mu(k\,\mu)^k\big]}{\tau^2\big[(\tau-s-\lambda)(\tau+k\,\mu)^k+\lambda(k\,\mu)^k\big]}.\\
   \end{aligned}
\end{equation}
From (\ref{S3_f_Er}) we see that the numerator of $f(\tau;s)$ has a double zero at $\tau=0$, so $f$ is analytic at $\tau=0$ and the poles $\tau=\tau(s)$ satisfy
\begin{equation}\label{S3_tau_pole}
(\tau-s-\lambda)(\tau+k\,\mu)^k+\lambda(k\,\mu)^k=0.
\end{equation}

Consider first the limit $x,t\rightarrow\infty$ with $1<t/x<\infty$. We define $F$ by
\begin{equation}\label{S3_F}
F(s,x)=\frac{1}{2\pi i}\int_{Br_\tau}f(\tau;s)\,e^{\tau \,x}\,d\tau.
\end{equation}
From (\ref{S3_tau_pole}) we have
\begin{equation}\label{S3_tau_pole2}
\lambda\left(\frac{k\,\mu}{\tau+k\,\mu}\right)^k=s-\tau+\lambda.
\end{equation}
By applying the residue theorem in (\ref{S3_F}) and using (\ref{S3_tau_pole2}), we have
\begin{equation}\label{S3_F_asym}
F(s,x)\sim R(\tau(s),s)\,e^{\tau(s)\,x},
\end{equation}
where
\begin{equation}\label{S3_R}
R(\tau(s),s)=\frac{s^2\,(\tau(s)+k\,\mu)}{\tau^2(s)\,\big[(k+1)\,\tau(s)-k\,s+k\,(\mu-\lambda)\big]}
\end{equation}
and $\tau=\tau(s)$ is the largest real root of (\ref{S3_tau_pole}). The other poles of $f$ lead to exponentially smaller terms and (\ref{S3_F_asym}) holds when $x\rightarrow\infty$ and $\Re(\tau)>0$ on the vertical contour $Br_\tau$.
We next define
\begin{equation}\label{S3_varphi}
\varphi(s)=\varphi\Big(s;\frac{t}{x}\Big)=s\,\frac{t}{x}-\tau(s),
\end{equation}
and then by (\ref{ptx_main}) we have
\begin{equation}\label{S3_ptx_1pole}
p(t|x)\sim\frac{1-\rho}{2\pi i}\int_{Br_s}e^{x\,\varphi(s)}\frac{1}{R(\tau(s),s)}\,ds.
\end{equation}
The integrand in (\ref{S3_ptx_1pole}) has saddle points where $\varphi'(s)=0$ so that there is a saddle point along the real axis at $\tilde{s}_*$, where $\tilde{s}_*$ and $\tilde{\tau}_*=\tau(\tilde{s}_*)$ are the solutions of
\begin{equation}\label{S3_s*tau*}
\left\{
   \begin{aligned}
& \varphi'(s)=\frac{t}{x}-\tau'(s)=0, \\
   & (\tau-s-\lambda)\,(\tau+k\,\mu)^k+\lambda(k\,\mu)^k=0.\\
  \end{aligned}
  \right.
\end{equation}
From (\ref{S3_tau_pole}) we have
\begin{equation}\label{S3_tau'}
\tau'(s)=\frac{\tau(s)+k\,\mu}{(k+1)\,\tau(s)+k(\mu-s-\lambda)}.
\end{equation}
Solving the system (\ref{S3_s*tau*}) with the help of (\ref{S3_tau'}) leads to (\ref{cor2.1.2s*}) and (\ref{cor2.1.2tau*}). If we shift $Br_s$ in (\ref{S3_ptx_1pole}) to $Br_s'$, on which $\Re(s)=\tilde{s}_*$, and use the saddle point method (see, e.g., Wong \cite{WO}) with the steepest decent direction $\arg(s-\tilde{s}_*)=\pm\pi/2$, we get
\begin{equation}\label{S3_ptx_saddle}
p(t|x)\sim\frac{1-\rho}{\sqrt{2\,\pi\, x\,\varphi''(\tilde{s}_*)}\;R(\tau(\tilde{s}_*),\tilde{s}_*)}\,e^{x\,\varphi(\tilde{s}_*)}.
\end{equation}
But, from (\ref{S3_varphi}) and (\ref{S3_tau_pole}) we have
\begin{equation}\label{S3_varphi''}
\varphi''(s)=-\tau''(s)=\frac{k\,(k+1)\,(\tau(s)+k\,\mu)\,(s-\tau(s)+\lambda)}{\big[(k+1)\,\tau(s)+k\,(\mu-s-\lambda)\big]^3}.
\end{equation}
Using (\ref{S3_R}), (\ref{S3_varphi}) and (\ref{S3_varphi''}) in (\ref{S3_ptx_saddle}) leads to (\ref{cor2.1.2}).

Next we consider $x,t\rightarrow\infty$ but with $t/x\approx 1$.
The previous calculation is not valid since, from
(\ref{cor2.1.2s*}), the saddle point
$\tilde{s}_*\rightarrow\infty$. From (\ref{S3_tau_pole}),
$\tilde{\tau}_*$ has the following expansion as
$\tilde{s}_*\to\infty$:
\begin{equation*}\label{S3_tau*}
\tilde{\tau}_*=\tilde{s}_*+\lambda-\frac{\lambda(k\,\mu)^k}{\tilde{s}_*^k}+O\Big(\frac{1}{\tilde{s}_*^{k+1}}\Big).
\end{equation*}
We return to (\ref{S3_ptx_1pole}) and note that $R(\tau(\tilde{s}_*),\tilde{s}_*)\sim 1$ as $\tilde{s}_*\rightarrow\infty$. Then we approximate the integrand for $s$ large (more precisely we can scale $s=O(x^{1/k})$ with $x(t-x)^k=O(1)$) to get the integral representation in (\ref{cor2.1.1}). Expanding $\exp[\lambda(k\,\mu)^ks^{-k}x]$ as a geometric series, using the identity
\begin{equation*}\label{S3_identity_case1}
\frac{1}{2\pi i}\int_{Br_\tau}e^{s\,(t-x)}\frac{1}{s^{km}}\,ds=  \left\{
   \begin{aligned}
    & \delta(t-x) & m=0, \\
    & \frac{(t-x)^{k\,m-1}}{(k\,m-1)!} & m\geq 1,  \\
   \end{aligned}
  \right.
\end{equation*}
and the generalized hypergeometric function
\begin{equation*}\label{S3_hypergeometric}
\sum_{m=0}^\infty\frac{z^m\,k^{k\,m}}{(m+1)!\,(k\,m+k-1)!}=\frac{1}{(k-1)!} \;_0F_k\big([\;];[1+\frac{1}{k},1+\frac{2}{k},...,2-\frac{1}{k},2];z\big),
\end{equation*}
we obtain the second expression in (\ref{cor2.1.1}).

Now we consider $x,t\rightarrow\infty$ but with $x/t$ small. From (\ref{cor2.1.2s*}) we let $x/t\rightarrow 0$, then the saddle point $\tilde{s}_*\rightarrow\tilde{s}_0$, which is given by (\ref{cor2.1.3s0}). Then from (\ref{cor2.1.2tau*}), $\tilde{\tau}_*$ has the following expansion
\begin{equation*}\label{S3_tau*expansion}
\tilde{\tau}_*=\tilde{\tau}_0\pm \tilde{\tau}_a\sqrt{s-\tilde{s}_0}+O(s-\tilde{s}_0),
\end{equation*}
where $\tilde{\tau}_0$ is given by (\ref{cor2.1.3tau0}) and
\begin{equation*}\label{S3_taua}
\tilde{\tau}_a=\sqrt{\frac{2\,k\,(\lambda\,\mu^k)^{\frac{1}{k+1}}}{k+1}}.
\end{equation*}
This means that on this scale, we must re-examine (\ref{S3_F}) as now two poles of the function $f(\tau;s)$ determine the asymptotics of the integral (\ref{S3_F}). We denote these two poles as $\tilde{\tau}_+=\tilde{\tau}_+(\tilde{s}_0)$ and $\tilde{\tau}_-=\tilde{\tau}_-(\tilde{s}_0)$. Then for $s\to\tilde{s}_0$
\begin{eqnarray}\label{S3_F_case3}
F(s,x) &\sim& R(\tilde{\tau}_+,\tilde{s}_0)\,e^{\tilde{\tau}_+\,x}+R(\tilde{\tau}_-,\tilde{s}_0)\,e^{\tilde{\tau}_-\,x}\nonumber\\
       &\sim& \frac{R_0}{\sqrt{s-\tilde{s}_0}}\,e^{\tilde{\tau}_0\,x}\,\big(e^{\tilde{\tau}_a\sqrt{s-\tilde{s}_0}\;x}-e^{-\tilde{\tau}_a\sqrt{s-\tilde{s}_0}\;x}\big),
\end{eqnarray}
where $R(\tilde{\tau}_\pm,\tilde{s}_0)\sim \pm\frac{R_0}{\sqrt{s-\tilde{s}_0}}$ and $R_0$ is given by
\begin{equation}\label{S3_R0}
R_0=\frac{1}{\sqrt{2\,k^3\,(k+1)}}\sqrt{\mu\,\rho^{\frac{1}{k+1}}}\,\bigg[\frac{\rho-(k+1)\,\rho^{\frac{1}{k+1}}+k}{1-\rho^{\frac{1}{k+1}}}\bigg]^2.
\end{equation}
We return to (\ref{ptx_main}), shift the contour $Br_s$ to $Br''_s$, which is slightly to the right of $\tilde{s}_0$, and use (\ref{S3_F_case3}), thus obtaining
\begin{eqnarray}\label{S3_ptx_case3}
p(t|x) &\sim& \frac{1-\rho}{2\pi i}\int_{Br''_s}\frac{e^{s\,t-\tilde{\tau}_0\,x}\,\sqrt{s-\tilde{s}_0}}{R_0\,\big(e^{\tilde{\tau}_a\sqrt{s-\tilde{s}_0}\;x}-e^{-\tilde{\tau}_a\sqrt{s-\tilde{s}_0}\;x}\big)}ds\nonumber\\
&=& \frac{1-\rho}{R_0}\,e^{-\tilde{\tau}_0\,x}\,\sum_{n=0}^\infty\frac{1}{2\pi i}\nonumber\\
&&\times\int_{Br''_s}\sqrt{s-\tilde{s}_0}\,\exp\Big[s\,t-(2n+1)\,\tilde{\tau}_a\,\sqrt{s-\tilde{s}_0}\;x\Big]\,ds.
\end{eqnarray}
Using the identity
\begin{equation*}\label{S3_identity_case3}
\begin{aligned}
&\frac{1}{2\pi i}\int_{Br''_s}\sqrt{s-\tilde{s}_0}\,\exp\Big[s\,t-(2n+1)\,\tilde{\tau}_a\,\sqrt{s-\tilde{s}_0}\;x\Big]ds\\
&=\frac{1}{4\sqrt{\pi}}\Big[\frac{(2n+1)^2\,\tilde{\tau}_a^2\,x^2}{t^{5/2}}-\frac{2}{t^{3/2}}\Big]\,\exp\Big[\tilde{s}_0\,t-\frac{(2n+1)^2\,\tilde{\tau}_a^2\,x^2}{4\,t}\Big]
\end{aligned}
\end{equation*}
in (\ref{S3_ptx_case3}), we obtain (\ref{cor2.1.3}).

Finally, we consider the case $x=O(1)$ and $t\rightarrow\infty$. Now all the $k+1$ poles of $f$ in (\ref{S3_f_Er}) contribute to the asymptotics. We denote these poles by $\tau_i=\tau_i(s),\;i=1,\cdots,k+1,$ which are the solutions of (\ref{S3_tau_pole}), and (\ref{S3_F}) evaluates to
\begin{equation}\label{S3_F_case4}
F(s,x)=\sum_{i=1}^{k+1}e^{\tau_i(s)\,x}\,R_i(s),
\end{equation}
where $R_i(s)=R(\tau_i(s),s)$ is the residue of $f(\tau;s)$ at $\tau=\tau_i(s)$.
Then (\ref{ptx_main}) becomes
\begin{equation}\label{S3_ptx_case4}
p(t|x)=\frac{1-\rho}{2\pi i}\int_{Br_s}e^{s\,t}\,\frac{1}{\sum_{i=1}^{k+1}e^{\tau_i(s)\,x}\,R_i(s)}ds.
\end{equation}
From (\ref{S3_ptx_case4}) we obtain (\ref{cor2.1.4}) by locating the pole $s=\tilde{s}_c(x)$ with the largest real part, which is the maximal solution of $F(s,x)=0$.

We can simplify (\ref{cor2.1.4}) for $x\to\infty$, which leads to an explicit expression in the matching region between Cases 3 and 4. As $x\to\infty$, we have $\tilde{s}_c\rightarrow\tilde{s}_0$ and the expansion
\begin{equation*}\label{S3_sc}
\tilde{s}_c=\tilde{s}_0+\frac{\tilde{A}}{x}+\frac{\tilde{B}}{x^2}+\frac{\tilde{C}}{x^3}+O\Big(\frac{1}{x^4}\Big).
\end{equation*}
We then expand (\ref{S3_tau_pole}) about $\tau=\tilde{\tau}_0$, to find that $\tilde{A}=0$, $\tilde{B}<0$ (which corresponds to $\tilde{s}_c<\tilde{s}_0$) and that $\tau$ has the following expansion:
\begin{eqnarray}\label{S3_tau_case4}
\tau &=&\tilde{\tau}_0\pm i\sqrt{\frac{2\,k\,\mu\,\rho^{\frac{1}{k+1}}}{k+1}\,|\tilde{B}|}\;\frac{1}{x}\nonumber\\
&&+\;\bigg[\frac{k+2}{3\,(k+1)}\,\tilde{B}\mp i\sqrt{\frac{k\,\mu\,\rho^{\frac{1}{k+1}}}{2\,(k+1)\,|\tilde{B}|}}\,\tilde{C}\bigg]\frac{1}{x^2}+O\Big(\frac{1}{x^3}\Big).
\end{eqnarray}
Then as in the analysis of Case 3, two poles of the function $f(\tau;s)$ dominate the expansion of $F(s,x)$. We denote these two conjugate poles as $\tau_1$ and $\tau_2=\overline{\tau_1}$. From (\ref{S3_R}) we have
\begin{eqnarray}\label{S3_R1}
R_1(s)&=&R(\tau_1(s),s)\nonumber\\ &=&-i\sqrt{\frac{\mu\,\rho^{\frac{1}{k+1}}}{2\,k^3\,(k+1)\,|\tilde{B}|}}\,\bigg[\frac{\rho-(k+1)\,\rho^{\frac{1}{k+1}}+k}{1-\rho^{\frac{1}{k+1}}}\bigg]^2 x\nonumber\\
&&+\frac{\Big[(k-4)\,\rho^{\frac{1}{k+1}}-(k+2)\Big]\Big[\rho-(k+1)\,\rho^{\frac{1}{k+1}}+k\Big]^2}{3\,k^2\,(k+1)\,(1-\rho^{\frac{1}{k+1}})^3}\nonumber\\
&&-i\sqrt{\frac{\mu\,\rho^{\frac{1}{k+1}}}{(2\,k\,|\tilde{B}|)^3}}\bigg[\frac{\rho-(k+1)\,\rho^{\frac{1}{k+1}}+k}{1-\rho^{\frac{1}{k+1}}}\bigg]^2+O\Big(\frac{1}{x}\Big),
\end{eqnarray}
and $R_2(s)= \overline{R_1(s)}$. Thus (\ref{S3_F_case4}) is approximately
\begin{eqnarray}\label{S3_F_case4_apprx}
F(s,x)&\sim& e^{\tau_1(s)\,x}\,R_1(s)+e^{\tau_2(s)\,x}\,R_2(s)\nonumber\\
      &\sim& 2\;\Re(e^{\tau_1(s)\,x}\,R_1(s)).
\end{eqnarray}
Expanding (\ref{S3_F_case4_apprx}) as $x\rightarrow\infty$ and noting that $\tilde{s}_c$ is a root of $F(s,x)=0$ in this limit, we obtain the expressions for $\tilde{B}$ and $\tilde{C}$ in (\ref{cor2.1_B}) and (\ref{cor2.1_C}). Then we can simplify (\ref{cor2.1.4}) to
\begin{equation}\label{S3_ptx_case4_apprx}
p(t|x)\sim\frac{(1-\rho)\,\exp\big(\tilde{s}_0\,t+\tilde{B}\,{t}/{x^2}+\tilde{C}\,{t}/{x^3}\big)}{\frac{d}{ds}\Big[2\;\Re(e^{\tau_1(s)\,x}\,R_1(s))\Big]\Big\arrowvert_{s=\tilde{s}_c}}.
\end{equation}
But,
\begin{equation}\label{S3_ds_case4}
\frac{d}{ds}\Big(e^{\tau_1(s)\,x}\,R_1(s)\Big)=\Big[R'_1(s)+R_1(s)\,x\,\tau'_1(s)\Big]\,e^{\tau_1(s)\,x}
\end{equation}
and
\begin{equation}\label{S3_R1'_case4}
R'_1(\tilde{s}_c)=O(x^2),
\end{equation}
\begin{equation}\label{S3_tau1'_case4}
\tau'_1(\tilde{s}_c)\sim -i\,\frac{k\,\mu\,\rho^{\frac{1}{k+1}}}{\pi\,(k+1)}\;x.
\end{equation}
Using (\ref{S3_R1'_case4}), (\ref{S3_tau1'_case4}) and (\ref{S3_R1}), (\ref{S3_ds_case4}) becomes
\begin{equation*}\label{S3_ds1@sc_case4}
\frac{d}{ds}\Big(e^{\tau_1(s)\,x}\,R_1(s)\Big)\Big\arrowvert_{s=\tilde{s}_c}\sim\frac{\mu^2\,\rho^{\frac{2}{k+1}}}{(k+1)^2\,\pi^2}\,\bigg[\frac{\rho-(k+1)\,\rho^{\frac{1}{k+1}}+k}{1-\rho^{\frac{1}{k+1}}}\bigg]^2\,x^3\,e^{\tilde{\tau}_0\,x},
\end{equation*}
which leads to (\ref{matching_cor2.1}).

We now consider general service densities $b(y)$. The basic scales in Theorem 2.1 are the same as those for the $E_k$ case in Corollary 2.1, but some of the definiting equations are more complicated, becoming transcendental rather than algebraic.  We consider the function (\ref{f(tau,s)}) and note that $\hat{b}(0)=1$ and $\hat{b}'(0)=-1/\mu$, so $\tau=0$ is not a pole of $f(\tau;s)$. The poles $\tau=\tau(s)$ of $f(\tau;s)$ now satisfy
\begin{equation}\label{S3_general tau}
\tau-s-\lambda\,(1-\hat{b}(\tau))=0.
\end{equation}

For the case $x,t\rightarrow\infty$, with $1<t/x<\infty$, the asymptotics are obtained analogously to the $E_k$ case and (\ref{S3_varphi}) still applies, but (\ref{S3_R}), (\ref{S3_tau'}) and (\ref{S3_varphi''}) must be replaced by
\begin{equation}\label{S3_general R}
R(\tau(s),s)=\frac{s^2}{\tau^2(1+\lambda\,\hat{b}'(\tau))},
\end{equation}
\begin{equation}\label{S3_general tau'}
\tau'(s)=\frac{1}{1+\lambda\,\hat{b}'(\tau)},
\end{equation}
and
\begin{equation*}\label{S3_general varphi''}
\varphi''(s)=\frac{\lambda\,\hat{b}''(\tau)}{\big(1+\lambda\,\hat{b}'(\tau)\big)^3}.
\end{equation*}
The first equation in (\ref{th2.1.2*}) follows from using $\varphi'(s)=0$, (\ref{S3_varphi}) and (\ref{S3_general tau'}).

For the case $x,t\rightarrow\infty$ with $t/x\approx 1$, we have to make some assumptions about the behavior of $b(y)$ as $y\to 0$. We assume that $b(y)\sim\alpha\, y^{\nu-1}\;(\alpha,\;\nu>0)$ for $y\to 0^+$. Then the Laplace transform of the service time satisfies $\hat{b}(\tau)\sim \alpha\,\Gamma(\nu)\,\tau^{-\nu}$ as $\tau\to\infty$. Then (\ref{th2.1.1}) is obtained in the same way as (\ref{cor2.1.1}), although we cannot express the sum in (\ref{th2.1.1}) as a hypergeometric function for non-integer $\nu$.

Next we consider $x,t\rightarrow\infty$ but with $x/t$ small. From (\ref{th2.1.2*}), letting $t/x\rightarrow\infty$, we have $s_0=s_*(\infty)$ and $\tau_0=\tau_*(\infty)$, which are given by (\ref{th2.1.3*}). Then $\tau_*$ has the following expansion as $s\to s_0$:
\begin{equation*}\label{S3_general tau expansion}
\tau_*=\tau_0\pm\sqrt{\frac{2(s-s_0)}{\lambda\,\hat{b}''(\tau_0)}}+O(s-s_0).
\end{equation*}
Thus again two poles of $f(\tau;s)$ dominate the expansion of $F(s,x)$ and the calculation is similar to the Erlang case, with (\ref{S3_R0}) becoming
\begin{equation*}\label{S3_general R0}
R_0=\frac{s_0^2}{\tau_0^2\,\sqrt{2\,\lambda\,\hat{b}''(\tau_0)}}.
\end{equation*}

Finally, we consider the case $x=O(1)$ and $t\rightarrow\infty$. For the general service time distribution, all the singularities of the function $f(\tau;s)$ contribute to $F(s,x)$. Then (\ref{th2.1.4}) is obtained by using the residue theorem at the largest pole $s_c(x)$ of the integrand in (\ref{ptx_main}), which is the maximal real solution of $F(s,x)=0$.

In the asymptotic matching region between Cases 3 and 4, we let $x\rightarrow\infty$ and
\begin{equation*}\label{S3_general Sc}
s_c=s_0+\frac{A}{x}+\frac{B}{x^2}+\frac{C}{x^3}+O\Big(\frac{1}{x^4}\Big),
\end{equation*}
and expand (\ref{S3_general tau}) at $\tau=\tau_0$. We find that $A=0$, $B<0$ and two conjugate poles $\tau_1$ and $\tau_2$ of the function $f(\tau;s)$ dominate the behavior of $F(s,x)$. Analogously to (\ref{S3_tau_case4}) and (\ref{S3_R1}), $\tau_1$ and $R_1(s)$ have the following expansions:
\begin{equation*}\label{S3_general tau1_R1}
\tau_1 = \tau_0+i\sqrt{\frac{2\,|B|}{\lambda\,\hat{b}''(\tau_0)}}\;\frac{1}{x}+\bigg[\frac{|B|\,\hat{b}'''(\tau_0)}{3\,\lambda\,(\hat{b}''(\tau_0))^2}-\frac{i\,C}{\sqrt{2\,\lambda\,\hat{b}''(\tau_0)\,|B|}}\bigg]\,\frac{1}{x^2}+O(\frac{1}{x^3}),
\end{equation*}
\begin{eqnarray*}
R_1(s) &=& -\frac{i\,s_0^2}{\tau_0^2\,\sqrt{2\,\lambda\,\hat{b}''(\tau_0)\,|B|}}\;x-\frac{s_0^2\,\big[6\,\hat{b}''(\tau_0)+\tau_0\,\hat{b}'''(\tau_0)\big]}{3\,\lambda\,\tau_0^3\,\big(\hat{b}''(\tau_0)\big)^2}\\
&&+\;\frac{i\,s_0^2\,C}{2\,B\,\tau_0^2\,\sqrt{2\,\lambda\,\hat{b}''(\tau_0)\,|B|}}+O\Big(\frac{1}{x}\Big),
\end{eqnarray*}
and $\tau_2 = \overline{\tau_1}$, $R_2(s) = \overline{R_1(s)}$. The constants $B$ and $C$ are obtained by expanding (\ref{S3_F_case4_apprx}) as $x\to\infty$ using $F(s_c(x),x)=0$, and this leads to (\ref{th2.1_BC}). Using (\ref{S3_ptx_case4_apprx}) with $(\tilde{s}_0,\tilde{B},\tilde{C})$ replaced by $(s_0,B,C)$ and
\begin{equation*}\label{S3_general tau'1}
\tau_1'(s_c)\sim -i\frac{x}{\pi\,\lambda\,\hat{b}''(\tau_0)},
\end{equation*}
\begin{equation*}\label{S3_general ds}
\frac{d}{ds}\Big(e^{\tau_1(s)\,x}\,R_1(s)\Big)\Big\arrowvert_{s=s_c}\sim\bigg(\frac{s_0}{\pi\,\lambda\,\tau_0\,\hat{b}''(\tau_0)}\bigg)^2\,x^3\,e^{\tau_0\,x},
\end{equation*}
we obtain (\ref{matching_th2.1}).

\section{Brief derivations of the conditional sojourn time density for the case $\rho\approx 1$}
Now we consider the $M/G/1$-PS model with a traffic intensity that is close to one, and let $\rho=1-\epsilon$ with $0<\epsilon\ll 1$.

First, we consider $x=O(1)$ and $t=O(1)$. Using $\lambda=\mu+O(\epsilon)$ we obtain from (\ref{f(tau,s)})
\begin{equation*}\label{S4_f expansion}
f(\tau;s)=\frac{s\,\big[\tau-\mu\,(1-\hat{b}(\tau))\big]}{\tau^2\,\big[\tau-s-\mu\,(1-\hat{b}(\tau))\big]}+O(\epsilon).
\end{equation*}
This leads to (\ref{th2.2.1}). On this scale the solution does not simplify much, but there is little probability mass in heavy traffic on the time scale $t=O(1)$.

Next, we consider $x=O(1)$ but for large time scales $t=T/\epsilon=O(\epsilon^{-1})$. In (\ref{S3_F}) we replace $\rho$ as $1-\epsilon$ and scale $s$ as $\epsilon w$, and we have
\begin{eqnarray*}\label{S4_F expansion}
F(s,x) &=& \frac{1}{2\pi i}\int_{Br_\tau}e^{\tau\, x}\bigg[\frac{\tau+w}{\tau^2}\,\epsilon+\frac{w^2}{\tau^2\,\big[\tau-\mu\,(1-\hat{b}(\tau))\big]}\,\epsilon^2+O(\epsilon^3)\bigg]\,d\tau\\
&=& \epsilon\, (1+wx)+\epsilon^2\,w^2\,\frac{1}{2\pi i}\int_{Br_\tau}\frac{e^{\tau\, x}}{\tau^2\,\big[\tau-\mu\,(1-\hat{b}(\tau))\big]}\,d\tau+O(\epsilon^3).
\end{eqnarray*}
Then from (\ref{ptx_main}), we obtain
\begin{equation}\label{S4_ptx case 2}
 \begin{aligned}
&p(t|x) = \frac{\epsilon^2}{2\pi i}\int_{Br_w}\frac{e^{w\,T}}{F(\epsilon\, w,x)}\,dw\\
&\sim \frac{\epsilon}{2\pi i}\int_{Br_w}e^{w\,T}\bigg[\frac{1}{1+wx}-\frac{\epsilon\, w^2}{2\pi i(1+wx)^2}\int_{Br_\tau}\frac{e^{\tau\, x}}{\tau^2\big[\tau-\mu(1-\hat{b}(\tau))\big]}d\tau\bigg]dw\\
&= \frac{\epsilon}{x}e^{-T/x}-\epsilon^2\bigg[\frac{\delta(T)}{x^2}+\frac{(T-2x)e^{-T/x}}{x^4}\bigg]\frac{1}{2\pi i}\int_{Br_\tau}\frac{e^{\tau\, x}}{\tau^2\big[\tau-\mu(1-\hat{b}(\tau))\big]}d\tau.
 \end{aligned}
\end{equation}
The function
\begin{equation*}\label{S4_g}
g(\tau,x)=\frac{e^{\tau\, x}}{\tau^2\,\big[\tau-\mu\,(1-\hat{b}(\tau))\big]}
\end{equation*}
has a pole at $\tau=0$ of order 4. By the residue theorem we have
\begin{equation}\label{S4_int g}
\frac{1}{2\pi i}\int_{Br_\tau}g(\tau,x)\,d\tau=Res_{\tau=0}\big(g(\tau,x)\big)+\frac{1}{2\pi i}\int_{\mathcal{C_-}}g(\tau,x)\,d\tau.
\end{equation}
Here we shifted the contour $Br_\tau$ to $\mathcal{C_-}$, which can be taken as the imaginary axis in the $\tau$-plane, indented to the left of $\tau=0$. Then we define
\begin{equation}\label{S4_Q*}
Q_*(x)=Res_{\tau=0}\big(g(\tau,x)\big)=\frac{1}{3!}\,\lim_{\tau\to 0}\frac{d^3}{d\tau^3}\bigg[\frac{\tau^2\,e^{\tau\, x}}{\tau-\mu\,(1-\hat{b}(\tau))}\bigg],
\end{equation}
which leads to (\ref{th2.2.2_Q*}). Note that we assumed that all the moments of the service time are finite, which are given by (\ref{moment}). Expression (\ref{th2.2.2}) is obtained by using (\ref{S4_ptx case 2}), (\ref{S4_int g}) and (\ref{S4_Q*}). The term proportional to $\delta(T)$ in (\ref{S4_ptx case 2}) does not mean that there is actually mass at $T=0$, but rather corresponds to the small ($O(\epsilon)$) mass that exists in the shorter time scale $t$, where (\ref{th2.2.1}) applies.

Now consider $x=X/\epsilon=O(\epsilon^{-1})$ and $t=T/\epsilon=O(\epsilon^{-1})$ with $1<T/X<\infty$. By the same argument as in Section 3, the pole $\tau=\tau(s)$ of $f(\tau;s)$ with the largest real part satisfies (\ref{S3_general tau}). We replace $\lambda$ by $\mu\,(1-\epsilon)$ in (\ref{S3_general tau}), which yields
\begin{equation}\label{S4_tau}
\tau-s-\mu\,(1-\epsilon)\,(1-\hat{b}(\tau))=0,
\end{equation}
and then expand $\tau$ as $\tau=\tau_a+\tau_b\,\epsilon+O(\epsilon^2)$. Then $\tau_a=\tau_a(s)$ and $\tau_b=\tau_b(s)$ satisfy
\begin{equation*}\label{S4_tau a}
\tau_a-s-\mu\,(1-\hat{b}(\tau_a))=0
\end{equation*}
and
\begin{equation*}\label{S4_tau b}
\tau_b=\frac{\mu\,(\hat{b}(\tau_a)-1)}{1+\mu\,\hat{b}'(\tau_a)}.
\end{equation*}
In (\ref{S3_varphi}) we replace $t$ and $x$ by $T/\epsilon$ and $X/\epsilon$ respectively, to get
\begin{eqnarray}\label{S4_varphi}
\varphi(s)=\varphi\Big(s;\frac{T}{X}\Big)&=&s\,\frac{T}{X}-\tau(s)\nonumber\\
&=&s\,\frac{T}{X}-\tau_a(s)-\tau_b(s)\,\epsilon+O(\epsilon^2)\nonumber\\
&=&\varphi_0(s)+\varphi_1(s)\,\epsilon+O(\epsilon^2)
\end{eqnarray}
and rewrite (\ref{S3_ptx_1pole}) as
\begin{equation}\label{S4_ptx 1pole}
p(t|x)\sim\frac{\epsilon}{2\pi i}\int_{Br_s}\frac{\exp\big[X\varphi_0(s)/\epsilon+X\varphi_1(s)\big]}{R(\tau(s),s)}\,ds.
\end{equation}
Here $R$ is as in (\ref{S3_general R}), with $\lambda$ replaced by $\mu$. Then the integrand in (\ref{S4_ptx 1pole}) has a saddle point where $\varphi'_0(s)=0$, which satisfies
\begin{equation}\label{S4_saddle}
  \left\{
   \begin{aligned}
    & \varphi_0'(s)=\frac{T}{X}-\tau'_a(s)=\frac{T}{X}-\frac{1}{1+\mu\,\hat{b}'(\tau_a)}=0,  \\
    & \tau_a-s-\mu\,(1-\hat{b}(\tau_a))=0.  \\
   \end{aligned}
  \right.
\end{equation}
We denote the solution of (\ref{S4_saddle}) as $\hat{\tau}_*=\hat{\tau}_*(T/X)$ and $\hat{s}_*=\hat{s}_*(T/X)$, which leads to (\ref{th2.2.4*}). Then by the standard saddle point method, (\ref{S4_ptx 1pole}) asymptotically evaluates to
\begin{equation}\label{S4_ptx case 4}
p(t|x)\sim\frac{\epsilon^{3/2}}{\sqrt{2\,\pi\, X\,\varphi''_0(\hat{s}_*)}\;R(\hat{\tau}_*,\hat{s}_*)}\exp\Big[X\varphi_0(\hat{s}_*)/\epsilon+X\varphi_1(\hat{s}_*)\Big].
\end{equation}
But by (\ref{S3_general R}) and (\ref{S4_varphi}), we have
\begin{equation}\label{S4_R case 4}
R(\hat{\tau}_*,\hat{s}_*)=\frac{\hat{s}_*^2}{\hat{\tau}_*^2\,\big[1+\mu\,\hat{b}'(\hat{\tau}_*)\big]}
\end{equation}
and
\begin{equation}\label{S4_varphi'' case 4}
\varphi''_0(\hat{s}_*)=-\tau_a''(\hat{s}_*)=\frac{\mu\,\hat{b}''(\hat{\tau}_*)}{\big[1+\mu\,\hat{b}'(\hat{\tau}_*)\big]^3}.
\end{equation}
Using (\ref{S4_R case 4}), (\ref{S4_varphi'' case 4}) and (\ref{S4_varphi}) in (\ref{S4_ptx case 4}), we obtain (\ref{th2.2.4}). We note that if $X\to\infty$ and $T\to\infty$ but $T/X=O(1)$, the approximation (\ref{th2.2.4}) remains valid.

For the case $X=O(1)$, $T=O(1)$ and $T-X\rightarrow 0^+$, we again assume that the service time density behaves as $b(y)\sim\alpha\, y^{\nu-1}\;(\alpha,\;\nu>0)$ for $y\to 0$. We note that the saddle point $\hat{s}_*\to\infty$ as $T/X\to 1$. Then from (\ref{S4_tau}), we find that $\tau$ has the following expansion for $s\to\infty$:
\begin{equation*}\label{S4_tau case 3}
\tau=s+\mu-\frac{\mu\,\alpha\,\Gamma(\nu)}{s^\nu}+O\Big(\frac{1}{s^{\nu+1}}\Big).
\end{equation*}
Following the same argument as in Section 3, we can easily obtain (\ref{th2.2.3}), once we scale $s$ as $S\epsilon^{1/\nu}$ and let $T_*=(T-X)\,\epsilon^{-1-1/\nu}=O(1)$.

Next, we consider $X=\sqrt{\epsilon}\,Z=O(\sqrt{\epsilon})$, $T=O(1)$. If we let $T/X\to\infty$ in (\ref{th2.2.4*}), it follows that $\hat{\tau}_*\to 0$ and the saddle point $\hat{s}_*\to 0$. By (\ref{S4_tau}) and scaling $s=O(\epsilon)=\epsilon w$, we find that $\tau$ has the following expansion:
\begin{equation*}\label{S4_tau case 5}
\tau\sim\pm\sqrt{\frac{2\,w}{\mu m_2}}\,\sqrt{\epsilon}+\frac{w\,m_3-3\,m_2}{3\,\mu m_2^2}\,\epsilon.
\end{equation*}
Now two poles of the function $f(\tau;s)$ dominate the behavior of $F(s,x)$. We approximate $F(s,x)$ by the sum of the residue at these two poles, where from (\ref{S3_general R}) we also have $R(\tau(s),s)\sim\pm\sqrt{\mu m_2w\epsilon/8}$. Then from (\ref{ptx_main}) we obtain
\begin{eqnarray*}\label{S4_ptx case 5}
p(t|x) &\sim& \frac{\sqrt{2}\,\epsilon^{3/2}}{\pi\, i\,\sqrt{\mu m_2}}\int_{Br_w}\frac{e^{w\,T}\exp\Big(-\frac{\sqrt{2\,w}\,Z}{\sqrt{\mu m_2}}\Big)}{\sqrt{w}\,\Big[1-\exp\Big(-\frac{2\,\sqrt{2\,w}\,Z}{\sqrt{\mu m_2}}\Big)\Big]}\,dw\\
&=& \frac{\sqrt{2}\,\epsilon^{3/2}}{\pi\, i\,\sqrt{\mu m_2}}\int_{Br_w}\frac{e^{w\,T}}{\sqrt{w}}\sum_{n=0}^\infty \exp\bigg[-\frac{(2n+1)\,\sqrt{2\,w}\,Z}{\sqrt{\mu m_2}}\bigg]\,dw,
\end{eqnarray*}
where the contour $Br_w$ is a vertical line in the $w$-plane slightly to the right of $w=0$.
Then (\ref{th2.2.5}) follows by using the identity
\begin{equation*}\label{S4_identity case 5}
\int_{Br_w}\frac{e^{w\,T}}{\sqrt{w}}\,\exp\bigg[-\frac{(2n+1)\,\sqrt{2\,w}\,Z}{\sqrt{\mu m_2}}\bigg]\,dw=\frac{1}{\sqrt{\pi\, T}}\exp\bigg[-\frac{(2n+1)^2\,Z^2}{2\,\mu m_2\,T}\bigg].
\end{equation*}
We note that by using the Poisson summation formula
\begin{equation*}\label{S4_poisson summation}
\sum_{n=-\infty}^\infty\psi(n)=\sum_{m=-\infty}^\infty\hat{\Psi}(2\pi m)=\sum_{m=-\infty}^\infty\int_{-\infty}^\infty e^{2\pi i ym}\psi(y)\,dy,
\end{equation*}
where $\hat{\Psi}$ is the Fourier transform of $\psi$, we can rewrite (\ref{th2.2.5}) as
\begin{equation}\label{S4_th2.2.5 rewrite}
p(t|x)\sim \frac{\epsilon^{3/2}}{Z}\bigg[1+2\sum_{n=1}^\infty(-1)^n\exp\Big(-\frac{n^2\,\pi^2\,\mu m_2\,T}{2\,Z^2}\Big)\bigg].
\end{equation}
From (\ref{S4_th2.2.5 rewrite}), we can easily verify that Cases 2 and 5 in Theorem 2.2 asymptotically match, in the intermediate limit where $x\rightarrow\infty$ and $Z\rightarrow 0$. Similarly, Cases 4 and 5 match in the limit where $X\rightarrow 0$ and $Z\rightarrow\infty$, which follows easily from (\ref{th2.2.5}).

Now we consider $X=O(1)$ and $T=\Theta/\epsilon=O(\epsilon^{-1})$ (thus $x=O(\epsilon^{-1})$ and $t=O(\epsilon^{-2})$). Similarly to the previous time scale, two poles, at $\tau_1=\tau_1(s)$ and $\tau_2=\tau_2(s)$, dominate the behavior of $F(s,x)$ and we have
\begin{eqnarray}\label{S4_th2.2.6}
p(t|x) &\sim& \frac{\epsilon}{2\pi i}\int_{Br_s}\frac{e^{s\,\Theta/\epsilon^2}}{Res_{\tau=\tau_1}(f(\tau;s)\,e^{\tau\, x})+Res_{\tau=\tau_2}(f(\tau;s)\,e^{\tau \,x})}\,ds\nonumber\\
&=& \frac{\epsilon}{2\pi i}\int_{Br_s}\frac{e^{s\,\Theta/\epsilon^2}}{R_1(s)\,e^{\tau_1(s)\,x}+R_2(s)\,e^{\tau_2(s)\,x}}\,ds.
\end{eqnarray}
We scale $s=O(\epsilon^{2})$ by setting $s=\big(\frac{\xi-1}{2\,\mu m_2}\big)\,\epsilon^2$, and then from (\ref{S4_tau}) and (\ref{S3_general R}) we have
\begin{equation*}\label{S4_case 6 tau_R}
\tau_{1,2}\sim\Big(\frac{-1\pm\sqrt{\xi}}{\mu m_2}\Big)\epsilon\;\;\textrm{and}\;\;R_{1,2}(s)=R(\tau_{1,2})\sim\pm\frac{(1\pm\sqrt{\xi})^2}{4\sqrt{\xi}}\epsilon,
\end{equation*}
which leads to (\ref{th2.2.6a}).

Furthermore, we expand the integrand in (\ref{th2.2.6a}) as a geometric series, and we have
\begin{equation}\label{S4_th2.2.6a geosum}
 \begin{aligned}
p&(t|x) \sim \frac{\epsilon^2}{\mu m_2\,\pi i}\exp\Big(\frac{X}{\mu m_2}-\frac{\Theta}{2\,\mu m_2}\Big)\\
  &\times \int_{Br_\xi}\frac{\exp\Big(\frac{\Theta\,\xi}{2\,\mu m_2}\Big)\sqrt{\xi}}{(1+\sqrt{\xi})^2}\sum_{n=0}^\infty\bigg(\frac{1-\sqrt{\xi}}{1+\sqrt{\xi}}\bigg)^{2n}\exp\Big[-\frac{(2n+1)X}{\mu m_2}\sqrt{\xi}\Big]d\xi.\\
 \end{aligned}
\end{equation}
Note that if we let $X\rightarrow\infty$, then the $n=0$ term in (\ref{S4_th2.2.6a geosum}) dominates, and we have
\begin{equation*}\label{S4_case 6 expansion}
 \begin{aligned}
&\frac{1}{2\pi i}\int_{Br_\xi}\frac{\sqrt{\xi}}{(1+\sqrt{\xi})^2}\exp\Big[\frac{\Theta}{2\mu m_2}\xi-\frac{X}{\mu m_2}\sqrt{\xi}\Big]d\xi\\
=&\bigg[\sqrt{\frac{2\mu m_2}{\pi\Theta}}+\sqrt{\frac{2\Theta}{\pi\mu m_2}}-\frac{2\mu m_2+X+\Theta}{\mu m_2}\,\mathrm{erfc}\left(\frac{X+\Theta}{\sqrt{2\mu m_2\Theta}}\right)\bigg]\exp\Big(-\frac{X^2}{2\mu m_2\Theta}\Big)\\
\sim&\; \sqrt{\frac{2\mu m_2}{\pi\Theta}}\exp\Big(-\frac{X^2}{2\mu m_2\Theta}\Big),\;\;\;X\rightarrow\infty.
 \end{aligned}
\end{equation*}
Here we used
$$\mathrm{erfc}(z)=1-\mathrm{erf}(z)=\frac{2}{\sqrt{\pi}}\int_z^\infty e^{-t^2}dt\sim\frac{1}{\sqrt{\pi}z}e^{-z^2},\textrm{ as }z\rightarrow\infty.$$
Then (\ref{S4_th2.2.6a geosum}) becomes, for $\Theta$ fixed and $X\rightarrow\infty$,
\begin{equation}\label{S4_case 5 matching}
p(t|x)\sim\frac{2^{3/2}\epsilon^2}{\sqrt{\pi\mu m_2\Theta}}\exp\Big[\frac{X}{\mu m_2}-\frac{\Theta}{2\mu m_2}-\frac{X^2}{2\mu m_2\Theta}\Big].
\end{equation}
When $X=O(\epsilon^{-1})$ and $T=O(\epsilon^{-1})$ but with $T/X=O(1)$, (\ref{th2.2.4}) remains valid, and letting $T/X\rightarrow\infty$ in (\ref{th2.2.4}) regains (\ref{S4_case 5 matching}). This again verifies that these two cases asymptotically match.

We return to (\ref{S4_th2.2.6a geosum}), and let $\sqrt{\xi}=z-1$, with which the integral becomes
\begin{eqnarray*}\label{S4_th2.2.6c}
  \begin{aligned}
2\sum_{n=0}^\infty & \exp\Big[\frac{\Theta}{2\mu m_2}+\frac{(2n+1)X}{\mu m_2}\Big]\\
& \times\int_\mathcal{C_+}\frac{(z-1)^2\,(z-2)^{2n}}{z^{2n+2}}\exp\Big[\frac{\Theta}{2\mu m_2}z^2-\frac{\Theta+(2n+1)X}{\mu m_2}z\Big]dz.
  \end{aligned}
\end{eqnarray*}
Here the contour $\mathcal{C_+}$ can be taken as the imaginary axis in the $z$-plane, indented to the right of $z=0$. Using the binomial expansion
$$(z-2)^{2n}=\sum_{j=0}^{2n}\frac{(-1)^j\,(2n)!}{j!\,(2n-j)!}\,2^{2n-j}\,z^j, $$
(\ref{S4_th2.2.6a geosum}) leads to
\begin{equation}\label{S4_th2.2.6b(1)}
 \begin{aligned}
p(t|x)\sim &\frac{\epsilon^2}{\mu m_2}\sum_{n=0}^\infty\sum_{j=0}^{2n}\frac{(2n)!}{j!\,(2n-j)!}\,(-2)^K\exp\Big[\frac{2(n+1)X}{\mu m_2}\Big]\\
&\times\; \frac{1}{2\pi i}\int_\mathcal{C_+}\frac{(z-1)^2}{z^K}\exp\Big(\frac{\Theta}{2\mu m_2}z^2-A_nz\Big)dz,
 \end{aligned}
\end{equation}
where
$$ A_n=\frac{\Theta+(2n+1)X}{2\mu m_2}$$
and
$$ K=2n-j+2.$$
We express the integral in (\ref{S4_th2.2.6b(1)}) in terms of parabolic cylinder functions, using
$$ \frac{1}{2\pi i}\int_{Br}z^\nu e^{-wz+z^2/2}dz=\frac{1}{\sqrt{2\pi}}D_\nu(w)e^{-w^2/4}, $$
thus obtaining
\begin{equation}\label{S4_th2.2.6b(2)}
 \begin{aligned}
p&(t|x)\sim \frac{\epsilon^2}{\sqrt{2\pi}\,\mu m_2} \sum_{n=0}^\infty \exp\Big[\frac{2(n+1)X}{\mu m_2}-\frac{z_n^2}{4}\Big]\sum_{j=0}^{2n}\frac{(-1)^j\,(2n)!}{j!\,(2n-j)!}\,(-2)^K\\
&\times\Big(\frac{\Theta}{\mu m_2}\Big)^{\frac{K-1}{2}}\left[\frac{\mu m_2}{\Theta}D_{2-K}(z_n)-2\sqrt{\frac{\mu m_2}{\Theta}}D_{1-K}(z_n)+D_{-K}(z_n)\right],
 \end{aligned}
\end{equation}
where $$z_n=\frac{\Theta+(2n+1)X}{\sqrt{\mu m_2\Theta}}.$$
Replacing $2n-j$ by $l$, (\ref{S4_th2.2.6b(2)}) leads to (\ref{th2.2.6b}).

If we let $X\rightarrow 0$ and $\Theta\rightarrow 0$ with $X/\sqrt{\Theta}$ (thus $z_n$) fixed, the term with $l=0$ in (\ref{th2.2.6c}) dominates and we have
\begin{eqnarray}\label{S4_th2.2.6b(3)}
p(t|x) &\sim& \frac{2^{3/2}\epsilon^2}{\sqrt{\pi\mu m_2\Theta}}\sum_{n=0}^\infty\exp\Big[\frac{2(n+1)X}{\mu m_2}-\frac{z_n^2}{4}\Big]D_0(z_n)\nonumber\\
       &\sim& \frac{2^{3/2}\epsilon^2}{\sqrt{\pi\mu m_2\Theta}}\sum_{n=0}^\infty\exp\Big[-\frac{(2n+1)^2X^2}{2\mu m_2\Theta}\Big].
\end{eqnarray}
Here we note that $X^2/\Theta=Z^2/T$ and  used the fact that $D_0(w)=e^{-w^2/4}$. Since (\ref{S4_th2.2.6b(3)}) is the same as (\ref{th2.2.5}), we have shown that Case 5 is really a special case of Case 6 in Theorem 2.2.

Alternately, we can treat the problem on the $(X,\Theta)$ scale by evaluating (\ref{S4_th2.2.6}) using the residues of the integrand at all the poles $s_p$, which satisfy
\begin{equation}\label{S4_th2.2.6c sp}
Res_{\tau=\tau_1(s_p)}(f(\tau;s)\,e^{\tau\, x})+Res_{\tau=\tau_2(s_p)}(f(\tau;s)\,e^{\tau \,x})= 0.
\end{equation}
We let $s_p=A_1\epsilon+B_1\epsilon^2+C_1\epsilon^3+O(\epsilon^4)$. Then from (\ref{S4_tau}) and (\ref{S3_general R}) we find that $A_1=0$ and we have the following expansions:
\begin{equation*}\label{S4_th2.2.6c tau12}
\tau_{1,2} \sim \Big(-\frac{1}{\mu m_2}\pm \tau_c\Big)\epsilon+(\tau_d\pm\tau_e)\epsilon^2,
\end{equation*}
where
\begin{equation*}
 \begin{aligned}
 & \tau_c=\frac{\sqrt{1+2\mu m_2B_1}}{\mu m_2},\\
 & \tau_d=\frac{2m_3-3\mu m_2^2+\mu m_2 m_3 B_1}{3\mu^2 m_2^3},\\
 & \tau_e=\frac{3\mu^2m_2^3C_1+3\mu m_2^2-2m_3+3\mu m_2(\mu m_2^2-m_3)B_1}{3\mu^2 m_2^3\sqrt{1+2\mu m_2B_1}},
 \end{aligned}
\end{equation*}
and
\begin{equation*}\label{S4_th2.2.6c R12}
R_{1,2}(s) \sim \Big(\frac{1}{2}\pm r_a\Big)\epsilon\pm r_b\,\epsilon^2,
\end{equation*}
where
\begin{equation*}
 \begin{aligned}
 & r_a=\frac{1+\mu m_2B_1}{2\sqrt{1+2\mu m_2B_1}},\\
 & r_b=\frac{\mu B_1\big[2\mu m_2^2C_1-(3\mu m_2^2-m_3)B_1\big]}{6(1+2\mu m_2B_1)^{3/2}}.\\
 \end{aligned}
\end{equation*}
By expanding the left-hand side of (\ref{S4_th2.2.6c sp}) about $\epsilon=0$, we obtain
\begin{equation}\label{S4_th2.2.6c O(1)}
e^{-2\tau_cX}=\frac{2r_a+1}{2r_a-1}
\end{equation}
and
\begin{equation}\label{S4_th2.2.6c O(e)}
r_b=2(r_a^2-1/4)\,\tau_e\,X.
\end{equation}
Setting $\sqrt{1+2\mu m_2B_1}=u+iv\;(v\geq 0)$ we find that all the roots of (\ref{S4_th2.2.6c O(1)}) are on the imaginary axis, and with $u=0$ (\ref{S4_th2.2.6c O(1)}) becomes (\ref{th2.2.6c_v}). Denoting the $n^{\textrm{th}}$ positive solution by $v_n=v_n(X)$, we obtain $B_1=B_1(v_n)$ by (\ref{th2.2.6c_Bvn}). By solving (\ref{S4_th2.2.6c O(e)}) for $C_1$, we have
\begin{equation*}\label{S4_th2.2.6c C1}
C_1=\frac{3\mu m_2(\mu m_2-m_3)XB_1^2+\big[(3\mu m_2^2-2m_3)X+\mu m_2(3\mu m_2^2-m_3)\big]B_1}{\mu^2m_2^3(1-XB_1)},
\end{equation*}
which leads to (\ref{th2.2.6c_Cvn}) with the help of (\ref{th2.2.6c_Bvn}). Note that $v_n=v_n(X)$ have the following asymptotic expansions:
\begin{equation}\label{S4_vn as X-inf}
v_n=\bigg[\frac{\mu m_2}{X}-\frac{2(\mu m_2)^2}{X^2}+\frac{4(\mu m_2)^3}{X^3}\bigg]\,n\,\pi+O\Big(\frac{1}{X^4}\Big),\;\;\;X\to\infty,
\end{equation}
and
\begin{equation*}\label{}
 \begin{aligned}
 &v_1=\frac{\sqrt{2\mu m_2}}{\sqrt{X}}-\frac{\sqrt{2}}{12\sqrt{\mu m_2}}\sqrt{X}+\frac{11\sqrt{2}}{1440(\mu m_2)^{3/2}}X^{3/2}+O(X^{5/2}),\;\;\;X\to 0,\\
 &v_n=\frac{(n-1)\pi\mu m_2}{X}+\frac{2}{(n-1)\pi}-\frac{4}{(n-1)^3\pi^3\mu m_2}X+O(X^2),\;\;\;n\geq 2,X\to 0.
 \end{aligned}
\end{equation*}
Now (\ref{S4_th2.2.6}) becomes
\begin{equation*}\label{S4_th2.2.6c ptx}
p(t|x)\sim \epsilon\sum_{n=1}^\infty\frac{e^{s_p(v_n)\,\Theta/\epsilon^2}}{\frac{d}{ds}\Big[e^{\tau_1(s)\,x}\,R_1(s)+e^{\tau_2(s)\,x}\,R_2(s)\Big]\Big\arrowvert_{s=s_p}}.
\end{equation*}
But
\begin{equation*}\label{}
\frac{d}{ds}\Big(e^{\tau_1(s)\,x}\,R_1(s)\Big)\Big\arrowvert_{s=s_p}=\Big[R'_1(s_p)+R_1(s_p)\,x\,\tau'_1(s_p)\Big]\,e^{\tau_1(s_p)\,x}
\end{equation*}
and
\begin{equation*}\label{}
R'_1(s_p)\sim \Big(\frac{1}{2}+i\frac{v_n^2-1}{4v_n}\Big)\,\epsilon,
\end{equation*}
\begin{equation*}\label{}
\tau'_1(s_p)\sim -i\frac{1}{v_n\,\epsilon},
\end{equation*}
\begin{equation*}\label{}
R_1'(s_p)\sim -i\frac{(1+v_n^2)\,\mu m_2}{4\,v_n^3\,\epsilon}.
\end{equation*}
Thus we obtain
\begin{equation}\label{S4_th2.2.6c d/ds1}
 \begin{aligned}
\frac{d}{ds}&\Big(e^{\tau_1(s)\,x}\,R_1(s)\Big)\Big\arrowvert_{s=s_p(v_n)} \\
&\sim \frac{1}{\epsilon}\bigg[\frac{v_n^2-1}{4v_n^2}X-i\frac{2v_nX+\mu m_2(v_n^2+1)}{4v_n^2}\bigg]\exp\bigg(\frac{-1+iv}{\mu m_2}X\bigg)
 \end{aligned}
\end{equation}
and
\begin{equation}\label{S4_th2.2.6c d/ds2}
\frac{d}{ds}\Big(e^{\tau_2(s)\,x}\,R_2(s)\Big)\Big\arrowvert_{s=s_p(v_n)}=\overline{\frac{d}{ds}\Big(e^{\tau_1(s)\,x}\,R_1(s)\Big)}\Big\arrowvert_{s=s_p(v_n)}.
\end{equation}
Using (\ref{S4_th2.2.6c d/ds1}) and (\ref{S4_th2.2.6c d/ds2}), we define $G$ by
\begin{equation*}\label{}
\frac{1}{G(v_n)}=2\;\Re\bigg\{\frac{1}{\epsilon}\Big[\frac{v_n^2-1}{4v_n^2}X-i\frac{2v_nX+\mu m_2(v_n^2+1)}{4v_n^2}\Big]\exp\Big(\frac{-1+iv}{\mu m_2}X\Big)\bigg\},
\end{equation*}
which leads to (\ref{th2.2.6c_G}), and then we obtain (\ref{th2.2.6c}) with $s_d(v_n)=s_p(v_n)/\epsilon^2$.

If we consider even larger time scales, with $\Theta\gg 1$ (thus
$t\gg \epsilon^{-2}$), then the largest pole $s_p(v_1)$ dominates.
Here $v_1=v_1(X)$ is the unique root in the interval $(0, \mu
m_2\pi/X)$ of (\ref{th2.2.6c_v}). This leads to (\ref{after
th2.2}). The expression (\ref{th2.2.6c}) with (\ref{th2.2.6c_Sd})
applies for time scales up to $\Theta=O(\epsilon^{-1})$
($t=O(\epsilon^{-3})$), but for even larger time scales we may
need further term in (\ref{th2.2.6c_Sd}), e.g., the
$O(\epsilon^{4})$ correction to $s_p$. We will discuss this more
in Section 5.

The $M/E_k/1$-PS results in Corollary 2.2 follows from Theorem 2.2 by using the $j^{\textrm{th}}$ moment
\begin{equation*}\label{}
m_j=\frac{(k+j-1)!}{k!\;k^{j-1}\,\mu^j}.
\end{equation*}

\section{Brief derivations of the unconditional sojourn time density}
The structure of the unconditional sojourn time density is highly dependent on the tail behavior of the service density. First we assume the service time density function behaves as (\ref{1.3}) or (\ref{b(y)_exp}). For $\rho$ fixed and less than one, the major contribution to the integral in (\ref{p(t)}) will come from the asymptotic matching region between the scales $x=O(1)$ and $x=O(\sqrt{t})$, with $t\to\infty$. In this region, the conditional sojourn time density is given asymptotically by (\ref{matching_th2.1}).

For $1\leq r< 2$, using (\ref{matching_th2.1}) and (\ref{b(y)_exp}) in (\ref{p(t)}), the unconditional sojourn time density behaves asymptotically as
\begin{equation}\label{S5_th2.3_pt1}
p(t)\sim \alpha_0\,e^{s_0\,t}\int_0^t x^{q-3}\exp\Big(-\tau_0\,x-N\,x^r+\frac{B\,t}{x^2}+\frac{C\,t}{x^3}\Big)dx,
\end{equation}
where
\begin{equation*}\label{S5_th2.3_alpha}
\alpha_0=\frac{(1-\rho)\,M\,\pi^2\,\lambda^2\,\tau_0^2\,[\hat{b}''(\tau_0)]^2}{2\,s_0^2}.
\end{equation*}
Scaling $x=\xi t^{\frac{1}{r+2}}=O(t^{\frac{1}{r+2}})$, (\ref{S5_th2.3_pt1}) becomes
\begin{equation}\label{S5_th2.3_pt2}
p(t)\sim \alpha_0\,t^{\frac{q-2}{r+2}}\,e^{s_0\,t}\int_{\delta t^{-\frac{1}{r+2}}}^{t^{\frac{r+1}{r+2}}}\xi^{q-3}\,e^{-\phi(\xi,t)}\,d\xi,
\end{equation}
where $\phi(\xi,t)$ is given by (\ref{th2.3.2_phi}). Here $\delta>0$ so as to avoid integration through $\xi=0$ in the case of $q<3$, and $\delta\gg 1$. By using the Laplace method with the major contribution coming from $\xi_*=\xi_*(t)$, which satisfies $\phi_{\xi}=0$ or (\ref{th2.3.2_phi'}), (\ref{S5_th2.3_pt2}) becomes
\begin{equation}\label{S5_th2.3_pt3}
p(t)\sim\frac{\sqrt{2\pi}\,\alpha_0}{\sqrt{\phi_{\xi\xi}(\xi_*,t)}}\,\xi_*^{q-3}\,t^{\frac{q-2}{r+2}}\,e^{s_0\,t-\phi(\xi_*,t)},
\end{equation}
where
\begin{eqnarray}\label{S5_th2.3_phi''}
\phi_{\xi\xi}(\xi,t) &=& \Big[r(r-1)N\xi^{r-2}-\frac{6B}{\xi^4}\Big]\,t^{\frac{r}{r+2}}-\frac{12C}{\xi^5}\,t^{\frac{r-1}{r+2}}\nonumber\\
&\sim& \Big[r(r-1)N\xi^{r-2}-\frac{6B}{\xi^4}\Big]\,t^{\frac{r}{r+2}}.
\end{eqnarray}

If $r=1$, then by (\ref{th2.3.2_phi'}), $\xi_*=\big(\frac{2|B|}{N+\tau_0}\big)^{1/3}$. Using (\ref{th2.3.2_phi}) and (\ref{S5_th2.3_phi''}) with $\xi_*$ and $r=1$ in (\ref{S5_th2.3_pt3}), we obtain (\ref{th2.3.1}).

If $1<r<2$, then by (\ref{th2.3.2_phi'}) the leading term $\xi_0$ in the asymptotic expansion of $\xi_*$ satisfies $Nr\xi^{r-1}+\frac{2B}{\xi^3}=0$, which leads to $\xi_0=\big(\frac{2|B|}{Nr}\big)^{\frac{1}{r+2}}$. Then we can rewrite (\ref{S5_th2.3_pt3}) as
\begin{equation*}
p(t)\sim\frac{\sqrt{2\pi}\,\alpha_0}{\sqrt{\big[r(r-1)N\xi_0^{r-2}-\frac{6B}{\xi_0^4}\big]t^{\frac{r}{r+2}}}}\,\xi_0^{q-3}\,t^{\frac{q-2}{r+2}}\,e^{s_0\,t-\phi(\xi_*,t)},
\end{equation*}
which leads to (\ref{th2.3.2}). We give three terms of asymptotic expansion for $\xi_*$ and $\phi(\xi_*,t)$ in Theorem 2.3, as $t\to\infty$. We note that the third terms in these expansions are different according as $1<r\leq 3/2$ or $r>3/2$.

If $r=2$, the above analysis is still valid but we need to include the additional factor $\exp(Dt/x^4)$ in (\ref{matching_th2.1}), and then in (\ref{S5_th2.3_pt1}). The constant $D$ is obtained by refining the approximation (\ref{matching_th2.1}) so that it applies for $x=O(t^{1/4})$. Thus (\ref{S5_th2.3_pt2}) and (\ref{S5_th2.3_pt3}) become
\begin{eqnarray*}
p(t) &\sim& \alpha_0\,t^{\frac{q-2}{4}}\,e^{s_0\,t}\int_{\delta t^{-1/4}}^{t^{3/4}}\xi^{q-3}\exp\Big(\frac{D}{\xi^4}+\phi(\xi,t)\Big)d\xi\\
     &\sim& \frac{\sqrt{2\pi}\,\alpha_0}{\sqrt{\phi_{\xi\xi}(\xi_*,t)}}\,\xi_*^{q-3}\,t^{\frac{q-2}{4}}\exp\Big(s_0\,t-\phi(\xi_*,t)+\frac{D}{\xi_*^4}\Big)\\
     &\sim& \frac{\sqrt{2\pi}\,\alpha_0}{\sqrt{\phi_{\xi\xi}(\xi_0,t)}}\,\xi_0^{q-3}\,t^{\frac{q-2}{4}}\exp\Big(s_0\,t-\phi(\xi_*,t)+\frac{D}{\xi_0^4}\Big).
\end{eqnarray*}
Here $\xi_0=(|B|/N)^{1/4}$ and
\begin{equation*}
\xi_*=\xi_0+\frac{3NC-\tau_0B}{8NB}\;t^{-\frac{1}{4}}+O(t^{-\frac{1}{2}}),
\end{equation*}
\begin{equation*}
\phi(\xi_*,t)=2\sqrt{N|B|}\;t^{\frac{1}{2}}-\frac{NC+\tau_0B}{N^{1/4}|B|^{3/4}}\;t^{\frac{1}{4}}-\frac{(3NC-\tau_0B)^2}{16NB^2}+O(t^{-\frac{1}{4}}),
\end{equation*}
\begin{equation*}
\phi_{\xi\xi}(\xi_0,t)\sim 8Nt^{1/2}.
\end{equation*}
Thus, after simplification, we obtain (\ref{th2.3.3}).

Now we consider the unconditional distribution in the heavy traffic limit, again assuming that the service time density function behaves as (\ref{b(y)_exp}).

For the time scale $t=T/\epsilon=O(\epsilon^{-1})$, we use (\ref{th2.2.2}) in (\ref{p(t)}), which leads to (\ref{th2.4.1}) after we integrate from $x=0$ to $x=\infty$. 

To compute the unconditional density $p(t)$ on the scale $t=\sigma_*/\epsilon^{r+2}=O(\epsilon^{-r-2})$ with $1\leq r<2$, we use (\ref{after th2.2}) and (\ref{b(y)_exp}) in (\ref{p(t)}) with $\Theta=\sigma_*/\epsilon^r$. Scaling $x=X/\epsilon=O(\epsilon^{-1})$, (\ref{p(t)}) becomes, since $\epsilon t\to\infty$,
\begin{equation}\label{S5_th2.4.2_pt1}
p(t)\sim \frac{M}{\epsilon^q}\int_0^{\infty}G(v_1)X^q\exp\Big[-\frac{1}{\epsilon^{r}}\psi(X,\sigma_*)\Big]dX,
\end{equation}
where
\begin{eqnarray*}\label{}
\psi(X,\sigma_*) &=& NX^r-B_1(v_1(X))\sigma_*-C_1(v_1(X))\sigma_*\epsilon\\
                 &=& \psi_0(X,\sigma_*)+\psi_1(X,\sigma_*)\epsilon.
\end{eqnarray*}
Hence (\ref{S5_th2.4.2_pt1}) is a Laplace type integral, and the major contribution will come from where $\psi$ is minimal, which should satisfy
\begin{equation}\label{S5_th2.4.2_phi0d}
\frac{\partial}{\partial X}\psi_0(X,\sigma_*)=NrX^{r-1}-\frac{d}{dX}B_1(v_1(X))\sigma_*=0.
\end{equation}
But from (\ref{th2.2.6c_Bvn}) and (\ref{th2.2.6c_v}), we have
\begin{equation}\label{S5_th2.4.2_B1'}
\frac{d}{dX}B_1(v_1(X))=\frac{v_1^2(X)\,[v_1^2(X)+1]}{\mu\, m_2\,X\,[v_1^2(X)+1]+2\mu^2\,m_2^2}.
\end{equation}
Using (\ref{S5_th2.4.2_B1'}) in (\ref{S5_th2.4.2_phi0d}), we obtain (\ref{th2.4.2_sigma*}). This defines $X=X(\sigma_*)$ implicitly. Denoting the right-hand side of (\ref{th2.4.2_sigma*}) as $\Omega(X)$, we can verify that $\Omega'(X)>0$, so that $\Omega(X)$ is a monotonically increasing function. As we discussed in Section 4, $v_1(X)\sim\mu m_2\pi/X$ as $X\to\infty$ and $v_1(X)\sim\sqrt{2\mu m_2/X}$ as $X\to 0^+$. Then $\Omega(X)\to\infty$ as $X\to\infty$ and $\Omega(X)\to 0$ as $X\to 0^+$. Hence there is a unique positive root of $\Omega(X)=\sigma_*$, which we denote by $\hat{X}=\hat{X}(\sigma_*)$. Then we use the standard Laplace method in (\ref{S5_th2.4.2_pt1}) to get, for $1\leq r<2$,
\begin{equation}\label{S5_th2.4.2_pt3}
p(t)\sim \frac{\sqrt{2\pi}\,M\,\epsilon^{\frac{r}{2}-q}\,\hat{X}^q\,G(v_1(\hat{X}))}{\sqrt{\frac{\partial^2}{\partial X^2}\psi_0(\hat{X},\sigma_*)}}\exp\big[-\psi(\hat{X},\sigma_*)\,\epsilon^{-r}\big],
\end{equation}
which leads to (\ref{th2.4.2}).

If $r=2$, to compute the unconditional sojourn time density on the time scale $t=\sigma/\epsilon^4=O(\epsilon^{-4})$, we need to include the $D_1(v_1)\epsilon^2=O(\epsilon^2)$ term in $s_d(v_1)$ in (\ref{after th2.2}). $D_1(v_1)$ is obtained in the same way that $B_1(v_1)$ and $C_1(v_1)$ are derived, which we discussed in Section 4. Analogously to (\ref{S5_th2.4.2_pt1}) and (\ref{S5_th2.4.2_pt3}), with $\sigma_*$ replaced by $\sigma$ and $r=2$, we have
\begin{eqnarray}\label{S5_th2.4.3_pt}
p(t) &\sim& \frac{M}{\epsilon^q}\int_0^{\infty}G(v_1)\,X^q\,e^{-\psi(X,\sigma_*)\,\epsilon^{-2}}\,e^{D_1(v_1)\,\sigma}\,dX\nonumber\\
&\sim& \frac{\sqrt{2\pi}\,M}{\sqrt{\frac{\partial^2}{\partial X^2}\psi_0(\tilde{X},\sigma)}}\,\epsilon^{1-q}\,\tilde{X}^q\,G(v_1(\tilde{X}))\nonumber\\
&\times&\exp\Big\{-\psi(\tilde{X},\sigma)\,\epsilon^{-2}+D_1(v_1(\tilde{X}))\,\sigma+\frac{\big[\frac{\partial}{\partial X}\psi_1(\tilde{X},\sigma)\big]^2}{2\frac{\partial^2}{\partial X^2}\psi_0(\tilde{X},\sigma)}\Big\},
\end{eqnarray}
where $\tilde{X}=\tilde{X}(\sigma)$ satisfies (\ref{th2.4.2_sigma*}) with $\sigma_*$ replaced by $\sigma$ and $r=2$. This leads to (\ref{th2.4.3}).

Next we assume that the service time density function $b(y)$ has
finite support for $0\leq y\leq A$ and behaves as
(\ref{b(y)_support}) near the maximum of its support. As
$t\rightarrow\infty$ with fixed $\rho<1$, we remove the condition
on $x$ by using the results in (\ref{th2.1.4}). The main
contribution comes from $x=A$, and we have
\begin{eqnarray}\label{S5_p(t)A1}
p(t) &\sim& \int_0^A \alpha_*\,(A-x)^{\nu_*-1}(1-\rho)\,J(x)\,e^{s_c(x)t}\,dx\nonumber\\
     &\sim& (1-\rho)\,\alpha_*\,J(A)\,e^{s_c(A)\,t}\int_{-\infty}^A(A-x)^{\nu_*-1}e^{s_c'(A)\,(x-A)t}dx.
\end{eqnarray}
Setting $x=A-u/(s_c'(A)\,t)$, (\ref{S5_p(t)A1}) becomes
\begin{equation*}\label{}
p(t)\sim \frac{(1-\rho)\,\alpha_*\,J(A)\,e^{s_c(A)\,t}}{[s_c'(A)t]^{\nu_*}}\int_0^\infty u^{\nu_*-1}e^{-u}du.
\end{equation*}
Using $\int_0^\infty u^{\nu_*-1}e^{-u}du=\Gamma(\nu_*)$, we obtain (\ref{th2.5.1}).

For the heavy traffic case, we remove the condition on $x$ by using (\ref{th2.2.2}), and on the large time scale $t=O(\epsilon^{-1})$ we obtain (\ref{th2.5.2}). For even larger times with $T\to\infty$, by using the Laplace method (with the main contribution from $x=A$), (\ref{th2.5.2}) becomes
\begin{equation*}\label{S5_p(t)A2}
p(t)\sim \frac{\epsilon\,\alpha_*\,e^{-T/A}}{A}\int_{-\infty}^A(A-x)^{\nu_*-1}\,e^{T\,(x-A)/{A^2}}dx.
\end{equation*}
This leads to (\ref{after_th2.5}).

\appendix
\section{{Appendix}}

We will give a brief derivation of the Laplace transform of the conditional sojourn time distribution with deterministic service density $b(y)=\delta(y-1/\mu)$. This was derived by Ott (see (5.16) in \cite{OT}) and more recently in \cite{EG}. However, these authors use arguments that are specific to the case $G=D$. Here we point out that these results also follow easily from the general $M/G/1$-PS model.

We rewrite (5.16) in \cite{OT} as
\begin{equation}\label{A1}
\mathbf{E}[e^{-s \mathbf{V}(1/\mu)}]=\frac{(1-\rho)(\lambda+s)^2\,e^{-\rho-s/\mu}}{s^2+\lambda\,\big[s+(1-\rho)(\lambda+s)\big]\,e^{-\rho-s/\mu}},
\end{equation}
where we replaced $z$ in \cite{OT} by 1. To prove (\ref{A1}), in view of that (\ref{Ee-sV2}) we need to prove that
\begin{equation}\label{A2}
\frac{1}{2\pi i}\int_{Br_\tau}e^{\tau/\mu}f(\tau;s)d\tau=\frac{s^2\,e^{\rho+s/\mu}+\lambda\,\big[s+(1-\rho)(\lambda+s)\big]}{(\lambda+s)^2}.
\end{equation}
But, by (\ref{f(tau,s)}),
\begin{equation*}\label{}
f(\tau;s)=\frac{(1-\rho)\tau+s}{\tau^2}+\frac{s^2}{\tau^2\big[\tau-s-\lambda(1-e^{-\tau/\mu})\big]},
\end{equation*}
and we have
\begin{equation*}\label{}
\frac{1}{2\pi i}\int_{Br_\tau}e^{\tau/\mu}\frac{(1-\rho)\tau+s}{\tau^2}d\tau=1-\rho+s/\mu.
\end{equation*}
Thus (\ref{A2}) is equivalent to proving the following identity:
\begin{equation}\label{A3}
 \begin{aligned}
\frac{1}{2\pi i} & \int_{Br_\tau}\frac{e^{\tau/\mu}}{\tau^2\,\big[\tau-s-\lambda\,(1-e^{-\tau/\mu})\big]}d\tau\\
& =\frac{e^{\rho+s/\mu}}{(\lambda+s)^2}+\frac{(1-\rho)\lambda^2+(2-\rho)\,s\,\lambda}{(\lambda+s)^2\,s^2}-\frac{s+\mu\,(1-\rho)}{\mu \,s^2}.
 \end{aligned}
\end{equation}
If we scale $\tau=\mu\, T$ and set $w=\rho+s/\mu$, (\ref{A3}) becomes
\begin{equation}\label{A4}
\frac{1}{2\pi i}\int_{Br_T}\frac{e^T}{\mu^2\,T^2\big[T-w+\rho\, e^{-T}\big]}\,dT=\frac{e^w-w-1}{\mu^2\,w^2}.
\end{equation}

We shift the contour $Br_T$ to the right so that $\Re(T)>\Re(w)$. Then upon expanding the integrand in (\ref{A4}) as a geometric series and multiplying (\ref{A4}) by $\mu^2$ we must show that
\begin{equation}\label{A5}
\sum_{L=0}^\infty\frac{1}{2\pi i}\int_{Br_T}\frac{(-1)^L\rho^L}{T^2\,(T-w)^{L+1}}\,e^{(1-L)\,T}\,dT=\frac{e^w-w-1}{w^2}.
\end{equation}
For $L\geq 1$ we can close the integration contour in the right half of the $T$-plane, and these integrals all evaluate to zero. For $L=0$ we close in the left half-plane, where there is a simple pole at $T=w$ and a double pole at $T=0$. Calculating the residues leads to (\ref{A5}), thus proving (\ref{A1}).

\end{document}